\documentclass[graybox]{svmult}

\usepackage{type1cm}        						% activate if the above 3 fonts are
\usepackage{makeidx}         						% allows index generation
\usepackage{graphicx}       						% standard LaTeX graphics tool
\usepackage{multicol}        						% used for the two-column index
\usepackage[bottom]{footmisc}						% places footnotes at page bottom

\usepackage{newtxtext}       % 
\usepackage[varvw]{newtxmath}       % selects Times Roman as basic font
\usepackage{enumitem}
\usepackage{booktabs}

\makeindex             % used for the subject index
\makeatletter
\def\toclevel@title{-1}
\def\toclevel@author{0}
\makeatother

\usepackage[ruled,vlined,linesnumbered,algo2e,algosection]{algorithm2e} % Before cleveref
\usepackage[colorlinks,citecolor=blue]{hyperref}	% hyperlinks
\usepackage[capitalise,nameinlink]{cleveref}

\crefname{algorithm}{Algorithm}{Algorithms}
\DontPrintSemicolon
\makeatletter
\renewenvironment{algomathdisplay}
 {\[}
 {\@endalgocfline\vspace{-\baselineskip}\]\hfill\strut\par}
\makeatother
\SetKw{KwOr}{or}
\SetKw{KwAnd}{and}

\newenvironment{minproblem}[2][]{%
	\begin{aligned}
		\min_{#1} \quad& #2
		\\
		\text{s.t.} \quad&
		\begin{aligned}[t]% the option t forces proper alignment
		}{%
		\end{aligned}
	\end{aligned}
}

\usepackage{mathtools}

\DeclarePairedDelimiter\norm{\lVert}{\rVert}

\DeclarePairedDelimiterX\innerprod[2]{\langle}{\rangle}{#1,#2}

\providecommand\given{\nonscript\;\delimsize|\nonscript\;\mathopen{}}
\DeclarePairedDelimiterX\set[1]{\{}{\}}{#1}

\newcommand\aeon{\;\text{a.e.\ on}\;}
\newcommand\adjoint{^\star}
\newcommand\LL{\mathcal{L}}
\newcommand\N{\mathbb{N}}
\newcommand\R{\mathbb{R}}
\newcommand{\embeds}{\hookrightarrow}
\newcommand{\Uad}{U_{\mathrm{ad}}}
\newcommand{\Wad}{W_{\mathrm{ad}}}
\newcommand{\dist}{\operatorname{dist}}
\newcommand{\cl}{\operatorname{cl}}
\newcommand{\conv}{\operatorname{conv}}
\newcommand{\cone}{\operatorname{cone}}

\newcommand{\blambda}{{\boldsymbol{\lambda}}}
\newcommand{\bx}{\mathbf x}
\newcommand{\bz}{\mathbf z}
\newcommand{\bd}{\mathbf d}

\newcommand*{\doi}[1]{\href{https://doi.org/\detokenize{#1}}{doi: \detokenize{#1}}}

\DeclareMathAlphabet{\mathpzc}{OT1}{pzc}{m}{it}
\newcommand\oo{\mathpzc{o}}
\newcommand\OO{\mathpzc{O}}

\begin{document}

\title*{Bilevel Optimal Control: Theory, Algorithms, and Applications}
\titlerunning{Bilevel Optimal Control: Theory, Algorithms, and Applications}
% Use \titlerunning{Short Title} for an abbreviated version of
% your contribution title if the original one is too long
\author{Stephan Dempe \and Markus Friedemann \and Felix Harder \and Patrick Mehlitz \and Gerd Wachsmuth}
\authorrunning{S. Dempe, F. Harder, M. Friedemann, P. Mehlitz, G. Wachsmuth}
% Use \authorrunning{Short Title} for an abbreviated version of
% your contribution title if the original one is too long
\institute{Stephan Dempe \at TU Bergakademie Freiberg, Faculty of Mathematics and Computer Science, Freiberg, Germany, \email{dempe@math.tu-freiberg.de}
\and Markus Friedemann \at BTU Cottbus-Senftenberg, Institute of Mathematics, Cottbus, Germany, \email{markus.friedemann@b-tu.de}
\and Felix Harder \at  \email{felix.harder.0@gmail.com}
\and Patrick Mehlitz \at University of Duisburg-Essen, Faculty of Mathematics, Essen, Germany, \email{patrick.mehlitz@uni-due.de}
\and Gerd Wachsmuth \at BTU Cottbus-Senftenberg, Institute of Mathematics, Cottbus, Germany, \email{wachsmuth@b-tu.de}%
}
%
% Use the package "url.sty" to avoid
% problems with special characters
% used in your e-mail or web address
%
\maketitle

%\abstract*{bla} % abstract only online

\abstract{%
	In this chapter, we are concerned with inverse optimal control problems, 
	i.e., optimization models which are used to identify parameters in optimal control problems 
	from given measurements.
	Here, we focus on linear-quadratic optimal control problems with
	control constraints where the reference control plays the role of the parameter and has to be reconstructed.
	First, it is shown that pointwise M-stationarity,  
	associated with the reformulation of the hierarchical model as a so-called
	mathematical problem with complementarity constraints (MPCC) in function spaces,
	provides a necessary optimality condition under some additional assumptions 
	on the data.
	Second, we review two recently developed algorithms 
	(an augmented Lagrangian method and a nonsmooth Newton method)
	for the computational identification of M-stationary points of finite-dimensional MPCCs.
	Finally, a numerical comparison of these methods, based on instances 
	of the appropriately discretized inverse optimal control problem of our interest,
	is provided.
}%				% abstract online and in printed book

\section{Introduction}

For decades, bilevel optimization problems, where two decision makers have to solve interdependent
optimization problems in predefined order, have been studied intensively from a theoretical, algorithmical,
and numerical point of view as they possess numerous underlying applications in data science, economy,
finance, machine learning, or natural sciences while being difficult due to inherent
nonconvexity, nondifferentiability, and irregularity.
For an introduction to bilevel optimization, we refer the interested reader 
to the monographs \cite{Dempe2002,DempeKalashnikovPerezValdesKalashnykova2015}.
A recent survey of applications, theory, and numerical methods addressing bilevel optimization can be found in \cite{Dempe2020}.
Typically, the two decision levels of a bilevel optimization problem are referred to 
as the upper- and lower-level problem, respectively.
Let us briefly recall the decision order in bilevel optimization.
First, the upper-level decision maker (the \emph{leader}) chooses from his individual feasible set a variable.
The latter is passed to the lower-level decision maker (the \emph{follower}) who now can solve his optimization
problem which is parametric in the leader's variable. Finally, the follower passes the obtained global solutions
back to the leader who now can evaluate his objective function which usually depends on the leader's and follower's variable.
As soon as the solution set of the follower is not a singleton, the overall bilevel optimization problem is ill-posed
as the follower may have to decide which of his multiple global solutions is handed over to the leader.
In order to overcome this issue, 
several approaches are known in the literature, see e.g.\ \cite{Zemkoho2016} for an overview.

Bilevel optimal control problems, which are bilevel optimization problems where at least one of the decision
makers has to solve an optimal control problem involving ordinary or partial differential equations,
are investigated e.g.\ in \cite{MehlitzWachsmuth2020}. 
Therein, the authors give a literature review, discuss the existence of optimal solutions, 
and comment on the derivation of necessary optimality conditions.
The latter topic is also studied in \cite{DempeHarderMehlitzWachsmuth2022} where, among others,
different strategies on how to derive stationarity conditions via regularization and relaxation approaches are investigated.
A particular subarea of bilevel optimal control is so-called inverse optimal control, 
where parameters in optimal control problems have to be reconstructed from given observations
(e.g., noisy state-control pairs of the underlying optimal control problem).
Typically, the upper-level problem of an inverse optimal control problem 
comprises a target-type term within the objective function 
modeling the desire of finding parameters 
such that the associated state-control pair, 
given implicitly via the lower-level optimal control problem, 
and the available observations are close to each other in a certain sense.
Several interesting applications of inverse optimal control, 
e.g.\ in the context of human locomotion,
see \cite{AlbrechtLeiboldUlbrich2012,AnbrechtPassenbergSobotkaPeerBussUlbrich2010,AlbrechtUlbrich2017,MombaurTruongLaumond2010},
visualize the importance of this model paradigm.
That is why inverse optimal control problems
have been studied theoretically throughout the last decade, so that a number of existence results,
optimality conditions, and solution algorithms is available in the literature, 
see e.g.\ \cite{DempeHarderMehlitzWachsmuth2019,DempeHarderMehlitzWachsmuth2022,FriedemannHarderWachsmuth2023,HarderWachsmuth2019,HatzSchloederBock2012,HollerKunischBarnard2018}.

Bilevel optimization problems possessing a convex lower-level problem 
(i.e., where the follower's problem is a convex optimization problem for each choice of the leader's variable)
are closely related to so-called mathematical problems with complementarity constraints (MPCCs for short),
see \cite{DempeDutta2012}, which are inherently irregular problems 
due to the highly combinatorial structure of their feasible set.
Fundamentals of finite-dimensional complementarity-constrained programming can be found 
in the classical monographs \cite{LuoPangRalph1996,OutrataKocvaraZowe1998},
and we refer the interested reader to the classical papers \cite{ScheelScholtes2000,Ye2005} 
where an overview of problem-tailored stationarity conditions and constraint qualifications for MPCCs is presented.
These concepts, at least partially, can be transferred to the infinite-dimensional setting,
see e.g.\ \cite{Mehlitz2017,MehlitzWachsmuth2016:1,Wachsmuth2014:1}.
A very popular notion of stationarity, which addresses MPCCs, is M-stationarity which is based on
the so-called Mordukhovich (or limiting) normal cone, see e.g. \cite{Mordukhovich2018}, to the complementarity set.
In finite dimensions, it is well known that M-stationarity provides a necessary optimality condition
in the presence of a very mild constraint qualification, see \cite{FlegelKanzow2006} for this classical result
and \cite{Harder2020} for a modern view.
The situation is pretty much different in infinite dimensions.
On the one hand, the Mordukhovich normal cone to complementarity sets in Lebesgue and Sobolev spaces
is often comparatively large as the limiting procedure involved in its construction tends to
annihilate information on the biactive set, see \cite{HarderWachsmuth2017,Mehlitz2017,MehlitzWachsmuth2016:2}.
On the other hand, working with a pointwise version of M-stationarity comes along with technical difficulties as
proof strategies known from the finite-dimensional situation cannot be transferred to the more general
setting without additional assumptions.
In fact, to the best of our knowledge, pointwise M-stationarity for MPCCs in Lebesgue and Sobolev spaces
has been considered merely in \cite{HarderWachsmuth2022,OutrataJarusekStara2011,Wachsmuth2014:2}.

In this chapter, after some brief discussion about the employed notation and fundamentals of complementarity-constrained
optimization in \cref{sec:preliminaries}, we first investigate a class of inverse optimal control problems with lower-level control constraints
whose local minimizers satisfy pointwise M-stationarity conditions (of a certain associated single-level MPCC reformulation) 
in the presence of additional assumptions on the problem data, see \cref{sec:IOC}. 
These results are based on the recent paper \cite{HarderWachsmuth2022}.
\Cref{sec:computing_M_stationary_points} presents two numerical methods which can be used to compute M-stationary points
of finite-dimensional MPCCs. We start with an augmented Lagrangian method which encapsulates the variational difficulty
of the complementarity set merely in the subproblem solver - a projected gradient type algorithm.
The associated theory can be found in \cite{GuoDeng2022,JiaKanzowMehlitzWachsmuth2023}.
Afterwards, a nonsmooth Newton method is presented which directly solves a nonsmooth (and even discontinuous) reformulation of
the M-stationarity system, see \cite{HarderMehlitzWachsmuth2021}.
For both algorithms, we first recapitulate the basic ideas and present some pseudocode before commenting on theoretical
convergence guarantees.
In \cref{sec:numerics_IOC}, we compare both algorithms by means of the (discretized) inverse optimal control problem
discussed in \cref{sec:IOC} whose minimizers are known to be M-stationary (in the finite- and infinite-dimensional framework).
We close with some concluding remarks in \cref{sec:conclusions}.

\section{Preliminaries}\label{sec:preliminaries}

\subsection{Notation}

For some positive integer $n\in\N$, we make use of $I^n:=\set{1,\ldots,n}$.
For $j\in I^n$, $\mathtt e\in\R^n$ and $\mathtt e_j\in\R^n$
represent the all-ones vector and the canonical unit vector
with a $1$ at entry $j$.
Furthermore, for reals $a_1,\ldots,a_n\in\R$, 
$[a_i]_{I^n}$ is the vector in $\R^n$ with components $a_1,\ldots,a_n$.
Whenever $a\in\R^n$ is a (row) vector and $I\subset I^n$ is an index set,
then $a_I$ is the vector resulting from $a$ by deleting those components
whose (row) index does not belong to $I$. 

Let $X$ be some (real) Banach space.
Its norm will be represented by $\norm{\cdot}_X\colon X\to[0,\infty)$.
The topological dual of $X$ will be denoted by $X\adjoint$.
Let $A\subset X$ be some arbitrary set and $x\in X$ be some point.
Then $\cl A$, $\conv A$, and $\cone A$ stand for 
the closure, the convex hull, and the conic hull of $A$, respectively. 
Furthermore, $\dist(x,A):=\inf_{y\in A}\norm{y-x}_X$ is the distance of $x$ to $A$,
and $\Pi(x,A):=\{y\in A\,|\,\dist(x,A)=\norm{y-x}_X\}$ denotes the projector of $x$ onto $A$.
Whenever $Y$ is yet another (real) Banach space, 
then $\LL(X,Y)$ is the Banach space of all continuous linear operators
which map from $X$ to $Y$. 
For $S\in\LL(X,Y)$, 
$S\adjoint\in\LL(Y\adjoint,X\adjoint)$ represents the adjoint of $S$.
If $F\colon X\to Y$ is (Fr\'{e}chet) differentiable at $\bar x\in X$,
then $F'(\bar x)\in\LL(X,Y)$ is the (Fr\'{e}chet) derivative of $F$ at $\bar x$.
In case $Y:=\R$, 
we make use of $\nabla F(\bar x):=F'(\bar x)\adjoint 1$ to represent the gradient of $F$ at $\bar x$.
Let us fix some (real) Hilbert space $H$.
The inner product in $H$ is denoted by $\innerprod{\cdot}{\cdot}_H\colon H\times H\to\R$.
We equip Hilbert spaces with the norm which is induced by the inner product in canonical way.

Throughout, we equip $\R^n$ with the Euclidean inner product and the associated Euclidean
norm which, for simplicity, is denoted by $\norm{\cdot}\colon\R^n\to[0,\infty)$.
For brevity, the operators $\max,\min\colon\R^n\times\R^n\to\R^n$ are interpreted in
componentwise fashion. 
Let $A\subset\R^n$ be a closed set and fix $\bar x\in A$.
Then
\[	
	N_A(\bar x)
	:=
	\left\{
		\eta\in\R^n\,\middle|\,
		\begin{aligned}
		&\exists\{x^k\}_{k\in\N},\{y^k\}_{k\in\N}\subset\R^n,\,\exists\{\alpha_k\}_{k\in\N}\subset(0,\infty)\colon
		\\
		&\quad
			x^k\to\bar x,\,\alpha_k(x^k-y^k)\to\eta,\,y^k\in\Pi(x^k,A)\,\forall k\in\N
		\end{aligned}
	\right\}
\]
is referred to as the limiting normal cone to $A$ at $\bar x$, 
see \cite{Mordukhovich2018} for a detailed introduction, examples, representations, and
calculus rules addressing this variational object. 

For some open bounded set $\Omega\subset\R^\ell$, $\ell\in\N$, 
which is equipped with the Lebesgue measure, and $p\in[1,\infty)$,
$L^p(\Omega)$ is the standard Lebesgue space of all (equivalence classes of)
measurable scalar functions on $\Omega$ whose $p$-th power is Lebesgue integrable.
We equip $L^p(\Omega)$ with the usual norm.
For $u\in L^p(\Omega)$, we make use of
$\set{u>0}:=\set{\omega\in\Omega\mid u(\omega)>0}$
and
$\set{u=0}:=\set{\omega\in\Omega\mid u(\omega)=0}$,
and we note that these sets are well defined up to sets of measure zero.
Similarly, sets with non-vanishing bounds are defined.
Let $\max\colon L^p(\Omega)\times L^p(\Omega)\to L^p(\Omega)$ be the
superposition operator associated with $\max\colon\R^2\to\R$.
We use $H^1(\Omega)$ for the classical Sobolev space of all function from $L^2(\Omega)$
possessing weak derivatives in $L^2(\Omega)$ with respect to all variables, 
and this space is equipped with the usual norm.
Further, $H^1_0(\Omega)$ denotes the closure of $C^\infty_\textup{c}(\Omega)$,
the set of all arbitrarily often continuously differentiable functions with compact support in $\Omega$, 
with respect to the $H^1$-norm, and this space is equipped with the conventional
$H^1$-seminorm, i.e.,
\begin{equation*}
	\forall w\in H^1_0(\Omega)\colon\quad
	\norm{w}_{H_0^1(\Omega)}
	:=
	\norm{\nabla w}_{L^2(\Omega)^\ell}.
\end{equation*}
Finally, $H^{-1}(\Omega):=H^1_0(\Omega)\adjoint$ is used.

\subsection{Finite-dimensional MPCCs}

Throughout the subsection, we consider a (finite-dimensional) mathematical program with complementarity constraints
given by
\begin{equation*}\label{eq:MPCC}\tag{MPCC}
	\begin{aligned}
		\min_{x\in\R^n} \quad& f(x)
		\\
		\text{s.t.}\quad&
		\begin{aligned}[t]
			g(x)	&\leq 	0,
			&
			h(x)	&=		0,
			\\
			G(x) 	&\geq	0,
			&
			H(x) 	&\geq 	0,
			&
			G(x)^\top H(x) &=	0,
		\end{aligned}
	\end{aligned}
\end{equation*}
where $f\colon\R^n\to\R$, $g\colon\R^n\to\R^r$, $h\colon\R^n\to\R^s$, and $G,H\colon\R^n\to\R^t$
are twice continuously differentiable functions.

We exploit the MPCC-tailored Lagrangian function $L\colon\R^n\times\R^r\times\R^s\times\R^t\times\R^t\to\R$
of \eqref{eq:MPCC} which is given by
\[
	L(x,\lambda,\eta,\mu,\nu)
	:=
	f(x)+\lambda^\top g(x)+\eta^\top h(x)+\mu^\top G(x)+\nu^\top H(x).
\]
For a feasible point $\bar x\in\R^n$ of \eqref{eq:MPCC},
we exploit the subsequently defined index sets:
\[
	\begin{aligned}
	&I_g&		&:=I_g(\bar x)&		&:=\set{i\in I^r\,|\,g_i(\bar x)=0},&\\
	&I_{+0}&	&:=I_{+0}(\bar x)&	&:=\set{i\in I^t\,|\,G_i(\bar x)>0\,\land\,H_i(\bar x)=0},&\\
	&I_{0+}&	&:=I_{0+}(\bar x)&	&:=\set{i\in I^t\,|\,G_i(\bar x)=0\,\land\,H_i(\bar x)>0},&\\
	&I_{00}&	&:=I_{00}(\bar x)&	&:=\set{i\in I^t\,|\,G_i(\bar x)=0\,\land\,H_i(\bar x)=0}.&
	\end{aligned}
\]
Let us recall that $\bar x$ is referred to as W-stationary (weakly stationary) whenever 
there are multipliers $\lambda\in\R^r$, $\eta\in\R^s$, and $\mu,\nu\in\R^t$ 
such that the conditions
\begin{subequations}\label{eq:MPCC_stat}
	\begin{align}
		\label{eq:MPCC_stat_x}
			&0=\nabla_xL(\bar x,\lambda,\eta,\mu,\nu),\\
		\label{eq:MPCC_stat_lambda}
			&\lambda_{I_g}\geq 0,\,\lambda_{I^r\setminus I_g}=0,\\
		\label{eq:MPCC_stat_W}
			&\mu_{I_{+0}}=0,\,\nu_{I_{0+}}=0
	\end{align}
\end{subequations}
are satisfied.
The point $\bar x$ is called C-stationary (Clarke-stationary)
if the multipliers additionally satisfy
\begin{equation}
	\label{eq:MPCC_stat_C}
	\forall i\in I_{00}\colon\quad \mu_i\nu_i\geq 0.
\end{equation}
Whenever \eqref{eq:MPCC_stat_C} can be strengthened to
\begin{equation}\label{eq:MPCC_stat_M}
	\forall i\in I_{00}\colon\quad (\mu_i<0\,\land\,\nu_i<0)\,\lor\,\mu_i\nu_i=0,
\end{equation}
then $\bar x$ is called M-stationary (Mordukhovich-stationary). 
The set of all multipliers, which solve the system of M-stationarity 
associated with $\bar x$, will be denoted by $\Lambda_\textup{M}(\bar x)$.
If some multipliers, which solve the system \eqref{eq:MPCC_stat},
additionally satisfy
\begin{equation}\label{eq:MPCC_stat_S}
	\forall i\in I_{00}\colon\quad \mu_i\leq 0 \,\land\, \nu_i\leq 0,
\end{equation}
then $\bar x$ is referred to as S-stationary (strongly stationary).
All these notions of stationarity can be found already in \cite{Ye2005}, and
their relations are obvious by definition.

Here, we will rely on yet another notion of stationarity taken from
\cite[Definition~2.1]{HarderWachsmuth2022}.
We say that $\bar x$ is A$_\forall$-stationary
if for each set $Q \subset I_{00}$,
there exist multipliers $\lambda\in\R^r$, $\eta\in\R^s$, and $\mu,\nu\in\R^t$ (depending on $Q$),
such that \eqref{eq:MPCC_stat} together with
\begin{equation}
	\label{eq:A_forall}
	\forall i \in I_{00}\setminus Q \colon\quad \mu_i\leq 0
	\qquad\land\qquad
	\forall i \in Q \colon\quad \nu_i\leq 0
\end{equation}
is satisfied.
The less restrictive concept of A-stationarity (Abadie-stationarity), 
where multipliers $\lambda\in\R^r$, $\eta\in\R^s$, and $\mu,\nu\in\R^t$
have to exist such that \eqref{eq:MPCC_stat} and \eqref{eq:A_forall} hold for
\emph{some} set $Q\subset I_{00}$, dates back to \cite{FlegelKanzow2003}.
It is worth mentioning that A$_\forall$-stationarity is equivalent to so-called
(linearized) B-stationarity (Bouligand-stationarity), see \cite[Proposition~3.1(b)]{Harder2021}.
Clearly, each S-stationary point is also A$_\forall$-stationary.
Furthermore, \cite[Theorem~3.4]{Harder2021} shows that each A$_\forall$-stationary point
is also M-stationary, see \cite[Section~3]{Harder2020} as well.

It is well known that the S-stationary points of \eqref{eq:MPCC} correspond to the KKT (Karush--Kuhn--Tucker)
points of this problem.
However, whenever $\bar x$ is a local minimizer of \eqref{eq:MPCC}, its S-stationarity 
can only be guaranteed in the presence of comparatively restrictive constraint
qualifications like MPCC-LICQ which demands that the gradients
\[
	\nabla g_i(\bar x)\,(i\in I_g),\;
	\nabla h_i(\bar x)\,(i\in I^s),\;
	\nabla G_i(\bar x)\,(i\in I_{0+}\cup I_{00}),\;
	\nabla H_i(\bar x)\,(i\in I_{+0}\cup I_{00})
\]
are linearly independent. In contrast, one can show that a local minimizer of \eqref{eq:MPCC}, where
the fairly mild constraint qualification MPCC-GCQ,
an MPCC-tailored variant of Guignard's constraint qualification, is valid, is already M-stationary,
see \cite{FlegelKanzow2006}.
We mention that MPCC-GCQ is satisfied
whenever the data functions $g$, $h$, $G$, and $H$ are affine.
An elementary and self-contained proof of this fact was recently given in
\cite{Harder2020}.
For later reference, we briefly sketch the outline of this proof.
Let $\bar x$ be a local minimizer of \eqref{eq:MPCC}
satisfying MPCC-GCQ.
By definition of this constraint qualification,
one immediately gets that $\bar d = 0$ is a locally optimal solution of
the linearized MPCC
\begin{equation*}\label{eq:linMPCC}\tag{linMPCC}
	\begin{aligned}
		\min_{d\in\R^n} \quad& f'(\bar x) d
		\\
		\text{s.t.}\quad&
		\begin{aligned}[t]
			\forall i \in I_g \colon
			g'_i(x) d	&\leq 	0,
			&
			\forall i \in I^s \colon
			h'_i(\bar x) d	&=		0,
			\\
			\forall i \in I_{0+} \colon
			G'_i(\bar x) d &= 0 ,
			&
			\forall i \in I_{+0} \colon
			H'_i(\bar x) d &= 0 ,
			\\
			\forall i \in I_{00} \colon
			G'_i(\bar x) d 	&\geq	0,
			&
			\forall i \in I_{00} \colon
			H'_i(\bar x) d 	&\geq 	0,
			\\
			&&
			\mathllap{
				\forall i \in I_{00} \colon
				(G'_i(\bar x) d) (H'_i(\bar x) d)
			} &=	0.
		\end{aligned}
	\end{aligned}
\end{equation*}
In case that the biactive set $I_{00}$
is not empty,
this problem still contains a complementarity condition.
Given any subset $Q \subset I_{00}$,
we can restrict the feasible set of \eqref{eq:linMPCC}
by replacing the constraints on $I_{00}$
by the stricter requirements
\begin{align*}
	\quad\forall i \in I_{00}\setminus Q \colon
	G'_i(\bar x)d &\geq 0,
	&
	\quad\forall i \in I_{00}\setminus Q \colon
	H'_i(\bar x)d &=    0,
	\\
	\quad\forall i \in Q\colon
	G_i(\bar x)d &=    0,
	&
	\quad\forall i \in Q\colon
	H'_i(\bar x)d &\geq 0.
\end{align*}
This results in a linear program with minimizer $\bar d = 0$.
It is straightforward to check that the associated optimality conditions
(which are always satisfied due to linearity)
are precisely \eqref{eq:MPCC_stat} together with \eqref{eq:A_forall}.
As $Q\subset I_{00}$ was chosen arbitrarily, this, in fact, shows that
$\bar x$ is A$_\forall$-stationary
(a slightly weaker version of this result can already be found in
\cite[Theorem~3.4 and associated discussions]{FlegelKanzow2005}).
This leads to a set of $2^{\# I_{00}}$ multipliers, where $\# I_{00}$
is the cardinality of the biactive set.
Finally,
one uses an elementary geometric argument, see \cite[Lemma~3.2]{Harder2020}, to prove
that there exists a convex combination of these multipliers
which additionally satisfies \eqref{eq:MPCC_stat_M}.
In fact, even stronger stationarity conditions can be distilled via this approach
when employing the Poincar\'{e}--Miranda theorem, see \cite[Section~3]{Harder2021}
for a discussion. 

The interest in C-stationarity is mainly due to its significance in numerical complementarity-constrained
optimization. Indeed, relaxation and penalty methods turn out to produce C-stationary points of
MPCCs, see e.g.\ \cite{HoheiselKanzowSchwartz2013,KanzowSchwartz2015,LeyfferLopezNocedal2006,RalphWright2004}.
Some recent contributions show that global convergence to M-stationary points can be achieved 
by certain sequential quadratic programming methods, see \cite{BenkoGfrerer2016},
multiplier-penalty methods, see \cite{GuoDeng2022,JiaKanzowMehlitzWachsmuth2023},
and (nonsmooth) Newton-type methods, see \cite{HarderMehlitzWachsmuth2021}.
The numerical detection of strongly stationary points is rather difficult as the set
\[
	M_\textup{S}
	:=
	\left\{
		(a,b,\mu,\nu)\in\R^4\,\middle|\,
		\begin{aligned}
			&0\leq a\perp b\geq 0,\,a\mu=0,\,b\nu=0,
			\\
			&a=b=0\,\Rightarrow\,\mu,\nu\leq 0
		\end{aligned}
	\right\},
\]
which is associated with S-stationarity, is not closed. We note that its closure
\begin{equation}\label{eq:M_stat_set}
	M_\textup{M}
	:=
	\left\{
		(a,b,\mu,\nu)\in\R^4\,\middle|\,
		\begin{aligned}
			&0\leq a\perp b\geq 0,\,a\mu=0,\,b\nu=0,
			\\
			&\mu\nu=0\,\lor\,(\mu<0\land\nu<0)
		\end{aligned}
	\right\}
\end{equation}
is characteristic for M-stationarity, i.e., 
the numerical search for M-stationary points is far more promising.

To end this section, we would like to comment on a certain second-order condition 
associated with an M-stationary point of \eqref{eq:MPCC}. 
Therefore, assume that $\bar x$ is an M-stationary point of \eqref{eq:MPCC}.
For $(\lambda,\eta,\mu,\nu)\in\Lambda_\textup{M}(\bar x)$, we define
\[
	S(\bar x,\lambda,\mu,\nu)
	:=
	\left\{
		d\in\R^n\,\middle|\,
		\begin{aligned}
			g'_i(\bar x)d 	&=0,\quad i\in I_g^+,\\
			h'(\bar x)d		&=0,\\
			G'_i(\bar x)d	&=0,\quad i\in I_{0+}\cup I_{00}^{\pm\R},\\
			H'_i(\bar x)d	&=0,\quad i\in I_{+0}\cup I_{00}^{\R\pm},\\
			(G'_i(\bar x)d)(H'_i(\bar x)d)	&=0,\quad i\in I_{00}^{00}
		\end{aligned}
	\right\},
\]
where we used the additional index sets
\[
	\begin{aligned}
	&I_g^+&				&:=I_g^+(\bar x,\lambda)&			&:=\set{i\in I_g\,|\,\lambda_i>0},&\\
	&I_{00}^{\pm\R}&	&:=I_{00}^{\pm\R}(\bar x,\mu,\nu)&	&:=\set{i\in I_{00}\,|\,\mu_i\neq 0},&\\
	&I_{00}^{\R\pm}&	&:=I_{00}^{\R\pm}(\bar x,\mu,\nu)&	&:=\set{i\in I_{00}\,|\,\nu_i\neq 0},&\\
	&I_{00}^{00}&		&:=I_{00}^{00}(\bar x,\mu,\nu)&		&:=\set{i\in I_{00}\,|\,\mu_i=\nu_i=0}.&
	\end{aligned}
\]
Observe that $S(\bar x,\lambda,\mu,\nu)$ is the union of finitely many subspaces of $\R^n$ and
serves as a generalization of the so-called critical subspace from nonlinear programming (NLP).
We say that MPCC-SSOC, the MPCC-tailored strong second-order condition, 
holds at $\bar x$ with respect to $(\lambda,\eta,\mu,\nu)$ whenever
\[
	\forall d \in S(\bar x,\lambda,\mu,\nu)\setminus\set{0}\colon\quad
	d^\top \nabla_{xx}L(\bar x,\lambda,\eta,\mu,\nu)d>0
\]
is valid. This corresponds to the validity of the classical SSOC, see e.g.\ \cite{Robinson1980},
along certain NLP branches of \eqref{eq:MPCC}. 
We note that MPCC-SSOC does not serve as a second-order sufficient optimality condition for \eqref{eq:MPCC}
as M-stationarity does not rule out the existence of descent directions.
Further discussion is provided in \cite[Section~2.2]{HarderMehlitzWachsmuth2021}.

\section{Inverse optimal control}
\label{sec:IOC}
We consider the framework of inverse optimal control,
i.e., the task of identifying parameters (via an optimization approach based on given measurements)
in a given optimal control problem.
Naturally, inverse optimal control problems possess a bilevel structure.
We closely follow the presentation from
\cite{HarderWachsmuth2022}.
We start with the lower-level problem
depending on a parameter $w \in L^2(\Omega)$.
It is given by
the parametrized optimal control problem
\begin{equation*}
	\label{eq:OC_w}
	\tag{OC$(w)$}
	\begin{minproblem}[u \in L^2(\Omega)]{
		\frac12 \norm{S u - y_{\textup{d}}}_Y^2 + \frac\alpha2 \norm{u - w}_{L^2(\Omega)}^2
	}
	u \ge u_{\textup{a}} \quad\aeon\Omega.
	\end{minproblem}
\end{equation*}
Here,
the parameter $w$ acts as a control offset or reference control.
The associated upper-level problem is given by
\begin{equation*}
	\label{eq:IOC}
	\tag{IOC}
	\begin{minproblem}[u \in L^2(\Omega), w \in H_0^1(\Omega)]
		{
			F(u, w)
			:=
			c(u)
			+ \frac12 \norm{w}_{H_0^1(\Omega)}^2
			+ \innerprod{\zeta}{w}_{L^2(\Omega)}
		}
		& \text{$u$ solves \eqref{eq:OC_w}}, \\
		& w \ge w_{\textup{a}} \quad\aeon\Omega.
	\end{minproblem}
\end{equation*}
One interpretation is that we are going to identify
the reference control $w$ in the lower-level problem
via measurements of the optimal control.
In this case, the objective $c\colon L^2(\Omega)\to\R$ would be of tracking type,
e.g.,
$c(u) := \frac12 \norm{u - u_{\textup{o}}}_{L^2(\Omega)}^2$
for some given observed control $u_{\textup{o}}\in L^2(\Omega)$.
Let us note that we cannot apply the theory from
\cite{DempeHarderMehlitzWachsmuth2019,FriedemannHarderWachsmuth2023}
for the derivation of stationarity conditions and solution
algorithms applying to \eqref{eq:IOC} as the upper-level variable
in this problem is chosen from an infinite-dimensional space. 
	Problems of related type have been considered, exemplary,
	in \cite{DempeHarderMehlitzWachsmuth2019,HatzSchloederBock2012,Mehlitz2017b,Ye1997}
	from a theoretical point of view.
	Some applications of inverse optimal control address 
	parameter reconstruction in optimal control problems appearing
	in the context of human locomotion,
	see e.g.\ \cite{AlbrechtLeiboldUlbrich2012,MombaurTruongLaumond2010},
	or scheduling,
	see e.g.\ \cite{PalagachevGerdts2017}.

We fix the assumptions concerning the data
in both problems.
The problem is posed on an open bounded set
$\Omega \subset \R^\ell$, $\ell \in \N$,
which is equipped with the Lebesgue measure.
Moreover,
we are given a (real) Hilbert space
$Y$, the control-to-state map $S \in \LL(L^2(\Omega), Y)$,
the desired state $y_{\textup{d}}\in Y$,
a regularization parameter $\alpha > 0$,
and
the lower-level lower bound
$u_{\textup{a}} \in L^2(\Omega)$.
In the upper-level problem,
the data satisfies
$\zeta \in L^2(\Omega)$
and
$w_{\textup{a}} \in H^1(\Omega)$ with $-\Delta w_{\textup{a}} \in L^2(\Omega)$
and $w_{\textup{a}} \le 0$ on $\partial\Omega$ in the sense $\max(w_{\textup{a}}, 0) \in H_0^1(\Omega)$.
The (partial) upper-level objective
$c \colon L^2(\Omega) \to \R$
is assumed to be continuously Fréchet differentiable
and bounded from below.
Further, we require some regularity of $S$.
Therefore, we assume the existence of $p > 2$
such that
$S\adjoint S \in \LL(L^2(\Omega), L^p(\Omega))$,
$S\adjoint y_{\textup{d}}, u_{\textup{a}} \in L^p(\Omega)$,
and $H_0^1(\Omega) \embeds L^p(\Omega)$.

We define the lower-level and upper-level feasible sets by
\begin{equation*}
	\Uad := \set{v \in L^2(\Omega) \given v \ge u_{\textup{a}} \text{ a.e.}}
	,\qquad
	\Wad := \set{v \in H_0^1(\Omega) \given v \ge w_{\textup{a}} \text{ a.e.}}
	,
\end{equation*}
respectively.
It is well known that
for fixed $w \in L^2(\Omega)$,
the unique solution $u \in \Uad$ of \eqref{eq:OC_w}
is characterized by the variational inequality
\begin{equation}
	\label{eq:OC_w_VI}
	\forall v \in \Uad\colon\quad
	\innerprod[\big]{S\adjoint(S u - y_{\textup{d}}) + \alpha (u - w)}{v - u}_{L^2(\Omega)}
	\ge
	0
	.
\end{equation}
The associated solution operator is continuous as a mapping from $L^2(\Omega)$ to $L^2(\Omega)$,
which can be seen by standard arguments.
Combining this with the compactness of the embedding $H_0^1(\Omega) \embeds L^2(\Omega)$,
see \cite[Theorem~1.34]{Troianiello1987},
the existence of solutions of \eqref{eq:IOC} can be shown.
For the proof of this result, we refer to \cite[Lemma~5.1]{HarderWachsmuth2022}.

\begin{lemma}
	\label{lem:existence_IOC}
	The inverse optimal control problem \eqref{eq:IOC}
	possesses a global solution.
\end{lemma}

Let us note that due to the inherent non-convexity of \eqref{eq:IOC},
one cannot expect uniqueness of solutions.

Due to the simple structure of the lower-level problem,
we can replace it with the associated KKT system.
This results in the KKT reformulation
\begin{equation*}
	\label{eq:KKTR}
	\tag{KKTR}
	\begin{minproblem}[u,\xi\in L^2(\Omega),w\in H_0^1(\Omega)]{
			F(u,w)
		}
		& w\geq w_{\textup{a}}\quad\aeon\Omega,
		\\
		& S\adjoint(Su-y_{\textup{d}})+\alpha(u-w)-\xi =0,\\
		& 0 \leq u-u_{\textup{a}} \perp \xi \geq0
		\quad\aeon\Omega,
	\end{minproblem}
\end{equation*}
which is an MPCC.
It can be shown that the multiplier $\xi$
is uniquely determined and depends continuously on $w$.
Thus,
\eqref{eq:IOC}
and
the above MPCC are equivalent
with respect to their local minimizers,
see e.g.\ \cite[Section~4.2.1]{Mehlitz2017}.

\subsection{M-stationarity for abstract problems in Lebesgue spaces}
\label{subsec:M_in_Lebesgue}
When considering infinite-dimensional MPCCs in Lebesgue or Sobolev spaces,
see e.g.\ \cite{GuoYe2016,HarderWachsmuth2017,HerzogMeyerWachsmuth2012,HintermuellerWegner2014,MehlitzWachsmuth2016:2,OutrataJarusekStara2011,Wachsmuth2014:2},
there are at least two possibilities
to generalize the M-stationarity conditions known from finite dimensions.
The first approach uses the fact that
the conditions \eqref{eq:MPCC_stat_W} and \eqref{eq:MPCC_stat_M}
encode that $(\mu, \nu)$ is an element of the limiting normal cone
to the complementarity set
\begin{equation}\label{eq:complementarity_set}
	\mathcal C
	:=
	\set{ (a,b) \in \R^t \times \R^t \given a \ge 0 \,\land\, b \ge 0 \,\land\, a^\top b = 0 }
	.
\end{equation}
This can be directly transferred to problems containing
the infinite-dimensional complementarity sets
\begin{equation*}
	\set{
		(u,\xi) \in L^2(\Omega) \times L^2(\Omega)
		\given
		u \ge 0 \text{ a.e.\ on $\Omega$} \,\land\,
		\xi \ge 0 \text{ a.e.\ on $\Omega$} \,\land\,
		\innerprod{u}{\xi}_{L^2(\Omega)} = 0
	}
\end{equation*}
or
\begin{equation*}
	\set*{
		(u,\xi) \in H_0^1(\Omega) \times H^{-1}(\Omega)
		\given
		\begin{aligned}
		&u \ge 0 \text{ a.e.\ on $\Omega$} \,\land\,
		\xi \ge 0 \text{ in the sense of $H^{-1}(\Omega)$}
		\\&\quad
		\,\land\,
		\xi(u) = 0
		\end{aligned}
	}
	.
\end{equation*}
However, it has been shown in
\cite{MehlitzWachsmuth2016:2}
that the associated limiting normal cone in the setting of Lebesgue spaces
is convex and does not contain any information on the biactive set,
i.e., one gets a system of W-stationarity.
In the setting of Sobolev spaces,
there is no precise characterization of the limiting normal cone available,
but it has been shown in \cite{HarderWachsmuth2017},
that (if the dimension $\ell$ of $\Omega$ is bigger than $1$)
the limiting normal cone is unreasonably large,
and it is very unlikely that it contains useful information
on the biactive set.

Another approach to define M-stationarity in infinite dimensions
is to require sign conditions as in 
\eqref{eq:MPCC_stat_W}, \eqref{eq:MPCC_stat_C}, \eqref{eq:MPCC_stat_M}, \eqref{eq:MPCC_stat_S},
or \eqref{eq:A_forall} in a pointwise (a.e.) sense.
In the remaining part of this subsection,
we will sketch the method from \cite{HarderWachsmuth2022} to prove such a pointwise M-stationarity
for an abstract class of MPCCs in Lebesgue spaces.
To the best of our knowledge, there are no other publications
providing pointwise M-stationarity for problems with complementarity constraints in Lebesgue spaces.
For a similar result in the one-dimensional setting with Sobolev spaces,
we are only aware of \cite{OutrataJarusekStara2011}, see also \cite{Wachsmuth2014:2}.

We start by
reproducing (a simplified version of) \cite[Theorem~3.1]{HarderWachsmuth2022}
without proof.
\begin{theorem}
	\label{thm:schinabeck}
	Let $\Omega_{00}\subset\Omega$ be measurable 
	and $\mathcal A$ be the Lebesgue-$\sigma$-algebra induced by $\Omega_{00}$.
	Further, let a family
	$\set{(\mu^Q, \nu^Q)}_{Q \in \mathcal{A}} \subset L^2(\Omega_{00}) \times L^2(\Omega_{00})$
	be given such that
	\begin{align}
		\label{eq:mu_nu_leq1}
		\mu^Q \leq0\quad\aeon\Omega_{00}\setminus Q
		\qquad
		\land
		\qquad
		\nu^Q \leq0\quad\aeon Q
	\end{align}
	for all $Q \in \mathcal{A}$.
	We assume the existence of a measurable function 
	$c_0\colon\Omega_{00}\to\R$
	and a constant $c_1>0$ with
		\begin{align*}
			\max(\mu^Q,\nu^Q) &\leq c_0
			\qquad\aeon\Omega,
			\quad\forall Q\in\mathcal{A},
			\\
			\norm{\mu^Q}_{L^2(\Omega_{00})}+\norm{\nu^Q}_{L^2(\Omega_{00})}
			&\leq c_1
			\qquad\forall Q\in\mathcal{A}.
		\end{align*}
	Then, there exists a point 
	\begin{equation*}
		(\bar\mu,\bar\nu) \in
		\operatorname{cl\,conv}\set{ (\mu^Q,\nu^Q)
		\given Q\in\mathcal{A} }\subset L^2(\Omega_{00})\times L^2(\Omega_{00})
	\end{equation*}
	which satisfies
	\begin{equation}
		\label{eq:mstat}
		(\bar\mu<0\,\land\,\bar\nu<0)\,\lor\,\bar\mu\bar\nu=0
		\qquad\aeon\Omega_{00}.
	\end{equation}
\end{theorem}

\Cref{thm:schinabeck} is the key result
to pass from a family of multipliers
satisfying the condition \eqref{eq:mu_nu_leq1},
which resembles the characteristic condition \eqref{eq:A_forall}
of A$_\forall$-stationarity in finite dimensions,
to a pair of multipliers satisfying
the pointwise M-stationary condition
\eqref{eq:mstat}.

Thus, it remains to prove the existence of a family of multipliers
satisfying \eqref{eq:mu_nu_leq1}.
In the finite-dimensional situation,
a similar assertion followed from the linearized problem
\eqref{eq:linMPCC}
and from
the optimality conditions
of tightened linear programs.
In the infinite-dimensional situation,
we need constraint qualifications
to derive necessary optimality conditions even for linear programs.
Surprisingly,
this leads to the effect
that the construction of multipliers satisfying \eqref{eq:mu_nu_leq1},
i.e., the verification of A$_\forall$-stationarity,
is much harder (and requires more restrictive assumptions)
than the passage from A$_\forall$-stationarity
to M-stationarity.

Putting everything together,
this leads to the following result.
For the proof, we refer to
\cite[Theorems~4.9 and 4.10]{HarderWachsmuth2022}.
\begin{theorem}
	\label{thm:Mstat_Lebesgue}
	Let $(0,0,0)\in L^2(\Omega)^3$ be a local solution
	of the linearized MPCC
	\begin{equation*}
		% \label{eq:mpcc2}
		% \tag{MPCC$_{\text{lin}}$}
		\begin{minproblem}[d_u , d_w,d_\xi\in L^2(\Omega)]{%
				\innerprod{F_u}{d_u}_{L^2(\Omega)} + \innerprod{F_w}{d_w}_{L^2(\Omega)} + \innerprod{F_\xi}{d_\xi}_{L^2(\Omega)}
			}
			A d_u -\kappa d_w - d_\xi &= 0,
			\\
			d_w&\geq 0 \quad\aeon\Omega_w,
			\\
			d_u&= 0\quad\aeon\Omega_{0+},
			\\
			d_\xi &= 0\quad\aeon\Omega_{+0},
			\\
			0 \leq d_u \perp d_\xi &\geq 0\quad\aeon\Omega_{00}
			.
		\end{minproblem}
	\end{equation*}
	Here,
	$A \colon L^2(\Omega)\to L^2(\Omega)$
	is a linear and continuous operator,
	$F_u, F_w,F_\xi\in L^{2}(\Omega)$ are fixed,
	and $\kappa>0$ is a constant.
	Additionally, the sets $\Omega_{0+},\Omega_{00},\Omega_{+0}\subset\Omega$
	are a measurable disjoint partition of $\Omega$,
	and $\Omega_w\subset\Omega$ is measurable.
	We further assume that the operator $A$ 
	enjoys one of the following properties.
	\begin{enumerate}[label=(\roman*)]
		\item
			\label{asm:operator_A:nonneg}
			The operator $A$ satisfies $A\adjoint v\geq0$ a.e.\ on $\Omega$
			for all $v\in L^2(\Omega)$ with $v\geq0$ a.e.\ on $\Omega$,
			and $\Omega_{+0}$ has measure zero.
		\item
			\label{asm:operator_A:average_with_scalars}
			The operator $A$ is given by
			$Av := C_1v + C_2\innerprod{1}{v}_{L^2(\Omega)}$ for all $v\in L^2(\Omega)$,
			where $C_1,C_2>0$ are fixed reals and $1\in L^2(\Omega)$ is the indicator
			function of $\Omega$, i.e., the function which is constantly $1$ on $\Omega$.
	\end{enumerate}
	Then, there exist multipliers
	$\bar p,\bar\lambda,\bar\mu,\bar\nu\in L^2(\Omega)$
	satisfying
	\[
			F_u + A\adjoint \bar p + \bar\mu = 0,
			\quad
			F_w - \kappa \bar p + \bar\lambda = 0,
			\quad
			F_\xi - \bar p + \bar\nu = 0,
	\]
	as well as the sign conditions
	\[
		\begin{aligned}
			\bar\lambda &\leq 0 \quad\aeon\Omega_w,
			&
			\bar\lambda &= 0 \quad\aeon\Omega\setminus\Omega_w,
			\\
			\bar\mu&=0 \quad\aeon\Omega_{+0},
			&
			\bar\nu&=0 \quad\aeon\Omega_{0+},
			\\&&
			\mathllap{(\bar\mu<0\,\land\,\bar\nu<0)\,\lor\,\bar\mu\bar\nu}&=0
			\quad\aeon\Omega_{00},
		\end{aligned}
	\]
	i.e., the pointwise M-stationarity conditions are valid.
\end{theorem}
As already said,
the restrictive assumptions on $A$
in \cref{thm:Mstat_Lebesgue}
are required
to prove validity of the KKT conditions
for linear programs
which are tightenings
of the linearized MPCC under consideration.
In \cite[Example~4.8]{HarderWachsmuth2022},
one can find an example of such a linear program
which fails to possess multipliers.

\subsection{M-stationarity for the inverse optimal control problem}

Now, we are going to explain
how the results from \cref{thm:Mstat_Lebesgue}
can be applied to the inverse optimal control problem
\eqref{eq:IOC}.
In a first step, we have to linearize the associated
equivalent complementarity-constrained problem \eqref{eq:KKTR}.
To this end, we need two important ingredients.

First, we have to check that the solution operator
of the lower-level problem is directionally differentiable.
This follows from
realizing that \eqref{eq:OC_w_VI}
can be rewritten as a projection
onto the feasible set $\Uad$
with respect to the non-standard inner product
associated with the coercive operator $\alpha \operatorname{id} + S\adjoint S$
on $L^2(\Omega)$.
Since the feasible set $\Uad$ is polyhedric,
the desired differentiability follows,
see \cite{Haraux1977,Mignot1976,Wachsmuth2016:2}.

Second, we need to show that local solutions $(\bar u, \bar w)\in L^2(\Omega)\times H^1_0(\Omega)$
of \eqref{eq:IOC}
enjoy the higher regularity $-\Delta \bar w \in L^2(\Omega)$.
This can be shown by a regularization approach,
see \cite[Theorem~5.4]{HarderWachsmuth2022}.

By combining these two steps,
we get that
$(0,0,0)\in L^2(\Omega)^3$ is a minimizer of
\begin{equation*}
	% \label{eq:linearized_OC}
	\begin{minproblem}
		[d_u, d_w, d_\xi \in L^2(\Omega)]
		{\innerprod{c'(\bar u)}{d_u}_{L^2(\Omega)} + \innerprod{-\Delta \bar w + \zeta}{d_w}_{L^2(\Omega)}}
		& (\alpha\operatorname{id} + S\adjoint S) d_u - \alpha d_w - d_\xi = 0, \\
		& d_w \ge 0 \aeon \set{\bar w = w_{\textup{a}}}, \\
		& d_u = 0 \aeon \set{\bar\xi > 0 }, \\
		& d_\xi = 0 \aeon \set{\bar u > u_{\textup{a}} }, \\
		& 0 \le d_u \perp d_\xi \ge 0 \aeon \set{ \bar u = u_{\textup{a}}} \cap \set{ \bar\xi = 0}
	\end{minproblem}
\end{equation*}
with
$\bar\xi := S\adjoint(S \bar u - y_{\textup{d}}) + \alpha (\bar u - \bar w)$,
whenever
$(\bar u, \bar w)\in L^2(\Omega)\times H^1_0(\Omega)$ is a local minimizer of \eqref{eq:IOC}.
Clearly, $\set{\bar u > u_\textup{a}},\set{\bar\xi>0},\set{\bar u=u_{\textup{a}}}\cap\set{\bar\xi=0}\subset\Omega$
provide a disjoint partition of $\Omega$ 
by definition of the complementarity constraint in \eqref{eq:KKTR}.

Putting everything together,
we obtain the following result,
see \cite[Theorem~5.9]{HarderWachsmuth2022},
	which shows that, in some specific situations,
	local minimizers of \eqref{eq:IOC} are
	(pointwise) M-stationary. 
	Potentially, this kind of stationarity can also be verified in a more
	general setting but, at the moment, it is not clear
	to us how such a result can be proven.
\begin{theorem}
	\label{thm:mstat}
	Let $(\bar u, \bar w) \in L^2(\Omega) \times H_0^1(\Omega)$ 
	be a local minimizer of \eqref{eq:IOC}
	and let $\bar\xi := S\adjoint(S \bar u - y_{\textup{d}}) + \alpha (\bar u - \bar w)$.
	We require validity of one of the following assumptions.
	\begin{enumerate}[label=(\roman*)]
		\item
			\label{thm:mstat:nonneg}
			The operator $S\adjoint S$
			satisfies $S\adjoint S v \ge 0$ a.e.\ on $\Omega$ for
			all $v\in L^2(\Omega)$ with $v \ge 0$ a.e.\ on $\Omega$,
			and $\bar u=u_{\textup{a}}$ holds a.e.\ on $\Omega$.
		\item
			\label{thm:mstat:average}
			It holds $Y=\R$ and
			$Sv = C\innerprod{1}{v}_{L^2(\Omega)}$ for all $v\in L^2(\Omega)$,
			where $C>0$ is a constant.
	\end{enumerate}
	Then, $(\bar u,\bar w,\bar\xi)$
	satisfies the pointwise M-stationary conditions of \eqref{eq:KKTR},
	i.e., there exist multipliers
	$\bar p,\bar\lambda,\bar\mu,\bar\nu\in L^2(\Omega)$
	satisfying
	\[
		c'(\bar u) + (\alpha\operatorname{id} + S\adjoint S) \bar p + \bar\mu = 0,
		\quad
		-\Delta \bar w +\zeta - \alpha \bar p + \bar\lambda = 0,
		\quad
		- \bar p + \bar\nu = 0,
	\]
	as well as the sign conditions
	\begin{align*}
		\bar\lambda &\leq 0 \quad\aeon\set{\bar w = w_{\textup{a}}},
		&
		\bar\lambda &= 0 \quad\aeon \set{\bar w > w_{\textup{a}}},
		\\
		\bar\mu&=0 \quad\aeon \set{\bar u > u_{\textup{a}}},
		&
		\bar\nu&=0 \quad\aeon \set{\bar \xi > 0},
		\\&&
		\mathllap{(\bar\mu<0\land\bar\nu<0)\lor\bar\mu\bar\nu}&=0
		\quad\aeon \mathrlap{\set{\bar u = u_{\textup{a}}} \cap \set{\bar \xi = 0}.}
	\end{align*}
\end{theorem}

In \cite[Section~5.4]{HarderWachsmuth2022}, it is been demonstrated 
by means of an example that, under the assumptions of \cref{thm:mstat}, 
local minimizers of \eqref{eq:IOC} are not necessarily pointwise S-stationary in
the sense that some multipliers which solve the pointwise M-stationarity 
system already satisfy
\[
	\bar\mu \leq 0
	\,\land\,
	\bar\nu \leq 0
	\quad
	\aeon \set{\bar u = u_{\textup{a}}} \cap \set{\bar \xi = 0}.
\]

\section{Computing M-stationary points of MPCCs}\label{sec:computing_M_stationary_points}

In this section, we present two algorithms which are tailored to find M-stationary points
of a given finite-dimensional MPCC. 
We start with a (safeguarded) augmented Lagrangian method from \cite{GuoDeng2022,JiaKanzowMehlitzWachsmuth2023}
which encapsulates the variationally difficult structure of MPCCs in the associated subproblems.
These, however, can be solved up to the necessary level of quality by suitable proximal gradient algorithms,
see e.g.\ \cite{DeMarchi2023,DeMarchiThemelis2022,JiaKanzowMehlitz2023,KanzowMehlitz2022} for recent studies.
Our second method, originating from \cite{HarderMehlitzWachsmuth2021}, 
applies a globalized semismooth Newton method to a nonsmooth reformulation of the M-stationarity system as
a discontinuous system of equations. 
	We note that both approaches are tailored to the finite-dimensional situation,
	see \cite{HarderWachsmuth2018,HarderWachsmuth2017,MehlitzWachsmuth2016:2} for a critical discussion of
	M-stationarity in function space optimization as well as some comments about the difference between
	M-stationarity in the narrower sense, defined via the so-called limiting normal cone, 
	see \cite{Mordukhovich2018}, to the complementarity set,
	and pointwise M-stationarity.

\subsection{A multiplier-penalty approach}\label{sec:ALM}

The principle idea behind the method of interest is to encapsulate the difficult variational structure 
of the complementarity constraints in \eqref{eq:MPCC} within an abstract geometric constraint set. 
This can be easily achieved with the aid of additional slack variables.
Then, we apply a (safeguarded) augmented Lagrangian scheme to the resulting problem where, in the appearing
subproblem, the augmented Lagrangian function is minimized subject to the geometric constraints.
In this regard, our solution approach is similar to \texttt{ALGENCAN} from \cite{AndreaniBirginMartinezSchuverdt2008}
which applies to NLPs with additional abstract convex geometric constraint sets 
(such that projections onto these geometric constraint sets are numerically available) .

\subsubsection{The algorithm}\label{sec:ALM_alg}

Let us consider the problem
\begin{equation*}\label{eq:slackMPCC}\tag{slackMPCC}
	\begin{aligned}
		\min_{x\in\R^n,\,z_G,z_H\in\R^t} \quad& f(x)
		\\
		\text{s.t.}\quad&
		\begin{aligned}[t]
			g(x)	&\leq 	0,
			&
			h(x)	&=		0,
			\\
			G(x)-z_G 	&=	0,
			&
			H(x)-z_H	&= 	0,
			\\
			(z_G,z_H)&\in \mathcal C,
			&
			&
		\end{aligned}
	\end{aligned}
\end{equation*}
where $\mathcal C$, defined in \eqref{eq:complementarity_set} is, up to a permutation of components,
the $t$-fold Cartesian product of the standard complementarity set in $\R^2$. 
As the values of the slack variables are uniquely determined for any $x\in\R^n$, \eqref{eq:MPCC} and \eqref{eq:slackMPCC}
are equivalent with respect to global and local minimizers (in the natural way).
For a compact notation, we introduce 
\begin{align*}
	\bx &:= (x,z_G,z_H)\in\R^n\times\R^t\times\R^t,\\
	\blambda&:=(\lambda,\eta,\mu,\nu)\in\R^r\times\R^s\times\R^t\times\R^t,\\
	\hat\blambda&:=(\hat\lambda,\hat\eta,\hat\mu,\hat\nu)\in\R^r\times\R^s\times\R^t\times\R^t,
\end{align*}
in order to denote the variables, multipliers, and surrogate multipliers throughout the subsection.
Similar meanings are behind $\bx^k$, $\blambda^k$, and $\hat\blambda^k$, where $k\in\N$ is an iteration index.
We define $L^\textup{sl}\colon(\R^n\times\R^t\times\R^t)\times(\R^r\times\R^s\times\R^t\times\R^t)\to\R$,
the Lagrangian-type function associated with \eqref{eq:slackMPCC}, by means of
\begin{align*}
	L^\textup{sl}(\bx,\blambda)
	:={}
	&
	f(x)+\lambda^\top g(x)+\eta^\top h(x)+\mu^\top (G(x)-z_G)+\nu^\top(H(x)-z_H)
	\\
	={}
	&
	L(x,\lambda,\eta,\mu,\nu)-\mu^\top z_G-\nu^\top z_H.
\end{align*}
For some penalty parameter $\rho>0$, the associated \emph{augmented} Lagrangian function 
$L^\textup{sl}_{\rho}\colon(\R^n\times\R^t\times\R^t)\times(\R^r\times\R^s\times\R^t\times\R^t)\to\R$
is then given by
\begin{align*}
	L^\textup{sl}_{\rho}(\bx,\blambda)
	:=
	f(x) + \frac{\rho}{2}
	\Bigl(
		&\norm{\max(g(x)+\lambda/\rho,0)}^2 + \norm{h(x)+\eta/\rho}^2
		\\
		&+\norm{G(x)-z_G+\mu/\rho}^2+\norm{H(x)-z_H+\nu/\rho}^2
	\Bigr).
\end{align*}
By continuous differentiability of all involved data functions,
$L^\textup{sl}_\rho$ is a continuously differentiable function.
In our multiplier-penalty scheme, we need to keep track of feasibility (with respect to inequality and equality constraints)
of the iterates. At the same time, we aim to monitor approximate validity of the complementarity-slackness condition
with respect to the inequality constraints.
Therefore, we introduce a function $V_\rho\colon(\R^n\times\R^t\times\R^t)\times(\R^r\times\R^s\times\R^t\times\R^t)\to\R$ given by
\[
	V_\rho(\bx,\blambda)
	:=
	\max(\norm{\max(g(x),-\lambda/\rho)},\norm{h(x)},\norm{G(x)-z_G},\norm{H(x)-z_H}).
\]
Let us emphasize that $V_\rho(\bx,\blambda)=0$ holds if and only if 
$\bx$ satisfies the inequality and equality constraints in \eqref{eq:slackMPCC}
while \eqref{eq:MPCC_stat_lambda} is valid.

The pseudocode of our augmented Lagrangian method can be found in \cref{alg:ALM}.
By construction, \cref{alg:ALM} is a so-called safeguarded multiplier-penalty method since
the surrogate multipliers $\hat\blambda^k$ remain bounded 
as they are chosen to be projections of the
actual multipliers $\blambda^k$ onto a (very large) box, see \cref{item:projection_of_multipliers}. 
This boundedness is essential in order to obtain a global convergence result for \cref{alg:ALM}.
Indeed, if $\hat\blambda^k$ is replaced by $\blambda^k$
everywhere in \cref{alg:ALM}, then a somewhat classical augmented Lagrangian method is recovered,
but a satisfying global convergence theory for this method is not likely to exist,
see \cite{KanzowSteck2017}.

\begin{algorithm2e}[htb]
	\SetAlgoLined
	\KwData{
		parameters $ \rho_0 > 0$, $\gamma > 1$, $q_\textup{alm} \in (0,1)$, $C>0$, 
		$\{\varepsilon_k\}_{k\in\N}\subset(0,\infty)$, $\tau_{\textup{alm}}>0$,
		starting points $x^0 \in \R^n $ and $\blambda^0\in\R^r\times\R^s\times\R^t\times\R^t$\; 
	}
	{
		Set $k := 0$ and $\bx^0:=(x^0,0,0)$\;
	}
	\While{$k=0$ \KwOr $V_{\rho_{k-1}}(\bx^k,\blambda^{k-1})>\tau_{\textup{alm}}$}
	{
		{
		Set
		\begin{algomathdisplay}
			\begin{aligned}
				\hat\lambda^k&:=\max(0,\min(\lambda^k,C\mathtt e)),&	
					\quad
					\hat\eta^k&:=\max(-C\mathtt e,\min(\eta^k,C\mathtt e)),&
					\\
				\hat\mu^k&:=\max(-C\mathtt e,\min(\mu^k,C\mathtt e)),&	
					\quad
					\hat\nu^k&:=\max(-C\mathtt e,\min(\nu^k,C\mathtt e))&
			\end{aligned}
		\end{algomathdisplay}
		\label{item:projection_of_multipliers}
		}
		{
		Compute a point $\bx^{k+1}\in\R^n\times\R^t\times\R^t$ such that	
    	\begin{algomathdisplay}
       		\dist(
       		-\nabla_{\bx}L^\textup{sl}_{\rho_k}(\bx^{k+1},\hat\blambda^k),
       		N_{\R^n\times\mathcal C}(\bx^{k+1})
       		\leq \varepsilon_{k+1}
    	\end{algomathdisplay}
    	\label{item:solve_subproblem_ALM}
    	}
    	{
 		Set 
		\begin{algomathdisplay}
			\begin{aligned}
        		\lambda^{k+1} &:=  \max(\rho_k g(x^{k+1})+\hat\lambda^k,0),&
        			\quad
        			\eta^{k+1} &:= \rho_k h(x^{k+1})+\hat\eta^k,&
        			\\
				\mu^{k+1} &:= \rho_k(G(x^{k+1})-z_G^{k+1})+\hat\mu^k,&
					\quad
					\nu^{k+1} &:= \rho_k(H(x^{k+1})-z_H^{k+1})+\hat\nu^k&
			\end{aligned}
		\end{algomathdisplay}
		 \label{item:update_multipliers}
    	}
    	\eIf{
    		$k=0$ \KwOr $V_{\rho_k}(\bx^{k+1},\hat\blambda^k)\leq q_\textup{alm}\,V_{\rho_{k-1}}(\bx^k,\hat\blambda^{k-1})$}{
    		Set $\rho_{k+1}:=\rho_k$\;
    	}
    	{
    		Set $\rho_{k+1}:=\gamma\rho_k$\;
    	}	
    	{
    		Set $k := k + 1$\;
    	}
	}
	\caption{Safeguarded augmented Lagrangian method for \eqref{eq:MPCC}.}
	\label{alg:ALM}
\end{algorithm2e}

\Cref{item:solve_subproblem_ALM} demands that the augmented Lagrangian subproblem
\[
	\begin{minproblem}[\bx\in\R^n\times\R^t\times\R^t]
		{
			L^\textup{sl}_{\rho_k}(\bx,\hat\blambda^k)
		}
		& \bx\in\R^n\times\mathcal C
	\end{minproblem}
\]
has to be solved up to a sufficiently small threshold $\varepsilon_{k+1}$ of approximate stationarity.
It has been shown in \cite[Section~3]{JiaKanzowMehlitzWachsmuth2023} that certain
projected gradient methods typically deliver such points, and although $\mathcal C$
is a set of highly combinatorial structure, projections onto $\mathcal C$ (which are not necessarily unique)
can be computed fast with the aid of ready-to-use formulas, see \cite[Section~5.1]{JiaKanzowMehlitzWachsmuth2023}.
Usually, the iterate $\bx^k$ computed in the prior iteration
is used as the initial guess in order to solve the subproblem.

One can easily check that
\cref{item:update_multipliers} implies
\[
	\nabla_{\bx}L^\textup{sl}_{\rho_k}(\bx^{k+1},\hat\blambda^k)
	=
	\nabla_{\bx}L(\bx^{k+1},\blambda^{k+1})
	.
\]
This explains the precise form of the multiplier update.

Finally, we note that the termination criterion of \cref{alg:ALM} ensures approximate
feasibility of the primal part associated with the final iterate in case where
\cref{alg:ALM} terminates.

The introduction of slack variables can be avoided
in the case that the functions $G$ and $H$ have simple structure in the sense
that it is easily possible to compute projections onto
\begin{equation}\label{eq:definition_of_D}
	\mathcal D
	:=
	\set{
		x \in \R^n
		\given
			G(x) 	\geq	0, \;
			H(x) 	\geq 	0, \;
			G(x)^\top H(x) =	0
	}
	.
\end{equation}
Then,
the slack variables $(z_G, z_H)$,
the associated multipliers $(\mu,\nu)$,
and the associated penalty terms in the augmented Lagrangian function
can be removed in the realization of
\cref{alg:ALM}.
Moreover,
$\R^n \times \mathcal C$
has to be replaced by $\mathcal D$
in \cref{item:solve_subproblem_ALM}
of \cref{alg:ALM}.
For further details,
the interested reader is referred to
\cite{JiaKanzowMehlitzWachsmuth2023}.

\subsubsection{Convergence guarantees}\label{sec:ALM_convergence}

Here, we briefly comment on the convergence behavior of \cref{alg:ALM} without giving any mathematical details.
Therefore, we assume that \cref{alg:ALM} produces an infinite sequence of iterates where $\varepsilon_k\downarrow 0$.
To start, we would like to note that each (primal) accumulation point of this sequence is 
a so-called \emph{asymptotically} stationary point of \eqref{eq:slackMPCC},
see \cite[Theorem~4.3]{JiaKanzowMehlitzWachsmuth2023}.
This means that there exist a primal sequence converging to the point of interest, 
a sequence of perturbation parameters tending to zero, 
and a (potentially unbounded) sequence of multipliers 
such that perturbed M-stationarity conditions are valid along the iterates,
see e.g.\ \cite{AndreaniHaeserSecchinSilva2019,Mehlitz2023,Ramos2021} for details about asymptotic stationarity
for MPCCs and \cite{Mehlitz2020} for a broader view on this topic.
We also note that, due to \cite[Lemma~3.1]{DeMarchiJiaKanzowMehlitz2023}, the asymptotically stationary
points of the complementarity-constrained problems 
\eqref{eq:MPCC} and \eqref{eq:slackMPCC} coincide (up to adding/removing certain components of the involved sequences).
It is well known that asymptotically stationary points of \eqref{eq:MPCC} are already M-stationary under
a very mild constraint qualification which is referred to as asymptotic regularity in the literature,
see \cite[Section~3.2]{Mehlitz2020}. More precisely, asymptotic regularity is valid whenever 
an MPCC-tailored variant of RCPLD (the so-called relaxed constant positive linear dependence constraint qualification)
is valid, see \cite[Lemma~2.7]{JiaKanzowMehlitzWachsmuth2023} or \cite[Lemma~3.12]{Mehlitz2023}
for this result and suitable references, or whenever the data functions $g$, $h$, $G$, and $H$ are affine,
see \cite[Theorem~5.3]{Mehlitz2020}. Roughly speaking, asymptotic regularity allows to take the limit in the
approximately satisfied M-stationarity conditions even in situations where the underlying sequence of
multipliers is unbounded. This idea dates back to the seminal paper \cite{AndreaniMartinezRamosSilva2016}.

Summing up the above arguments, \cref{alg:ALM} produces M-stationary points of \eqref{eq:MPCC} 
in the presence of a very mild qualification condition, 
see \cite[Theorem~3]{GuoDeng2022} and \cite[Corollary~4.4]{JiaKanzowMehlitzWachsmuth2023},  
even if the subproblems are merely solved up to approximate stationarity, 
and this result does not depend on the choice of the starting point.

\subsection{A nonsmooth Newton-type approach}\label{sec:SSN}

In this subsection, we demonstrate that the system of M-stationarity associated with \eqref{eq:MPCC} 
can be rewritten as a square system of nonsmooth equations. The latter can be solved with the aid
of a suitable Newton-type method based on generalized derivatives. 
Finally, a globalization strategy is discussed.

Classical nonsmooth Newton methods have been developed more than thirty years ago,
see e.g.\ \cite{Qi1993,QiSun1993}, and are often based on the concept of semismoothness which addresses
locally Lipschitzian functions, see \cite{Mifflin1977}.
However, as the nonsmooth system we are going to derive here is built from inherently discontinuous functions,
the concept of semismoothness does not apply. 
Instead, we will rely on so-called Newton differentiability which has been introduced in the context
of infinite-dimensional applications of mathematical optimization, see e.g.\ 
\cite{HintermuellerItoKunisch2002,Ulbrich2002}.

Recall that a given mapping $\Phi\colon\R^n\to\R^m$ is called Newton differentiable on some set $X\subset\R^n$
with Newton derivative $D\Phi\colon\R^n\to\R^{m\times n}$ whenever
\[
	\Phi(x+d)-\Phi(x)-D\Phi(x+d)d=\oo(\norm{d})\qquad \text{for }d\to 0
\]
holds for all $x\in X$. In case where even
\[
	\Phi(x+d)-\Phi(x)-D\Phi(x+d)d=\OO(\norm{d}^2)\qquad \text{for }d\to 0
\]
holds for all $x\in X$, $\Phi$ is said to be Newton differentiable of order $1$ on $X$ with Newton derivative $D\Phi$.
Finally, if for each $x\in X$, there is some $\varepsilon_x>0$ such that
\[
	\forall d\in\mathbb B_{\varepsilon_x}(0)\colon\quad
	\Phi(x+d)-\Phi(x)-D\Phi(x+d)d=0
\]
holds, then $\Phi$ is referred to as Newton differentiable of order $\infty$ on $X$ with Newton derivative $D\Phi$.
Clearly, whenever $\Phi$ is continuously differentiable, then it is also Newton differentiable on each subset of
$\R^n$ with Newton derivative $D\Phi:=\Phi'$. If $\Phi'$ is locally Lipschitzian, the order of Newton differentiability is $1$.
It has been shown in \cite[Lemma~2.11]{HarderMehlitzWachsmuth2021} that Newton differentiability enjoys a standard chain rule
which preserves the minimal order of Newton differentiability involved. 
Let us now assume that $m=n$ and that $\Phi$ is Newton differentiable on a set $X$, which comprises the roots of $\Phi$,
with Newton derivative $D\Phi\colon\R^n\to\R^{n\times n}$. Given a starting point $x^0\in\R^n$, we are interested in the
Newton-type scheme
\begin{equation}\label{eq:Newton_iteration}
	x^{k+1}:=x^k-D\Phi(x^k)^{-1}\Phi(x^k),\qquad k=0,1,\ldots,
\end{equation}
and \cite[Theorem~2.9]{HarderMehlitzWachsmuth2021} shows that whenever $\bar x\in\R^n$ is a root of $\Phi$ such that
the Newton derivative $D\Phi$ takes uniformly invertible values locally around $\bar x$, then each infinite sequence computed
by this scheme converges superlinearly to $\bar x$ if only $x^0$ is sufficiently close to $\bar x$.
Further, the convergence is quadratic if the order of Newton differentiability is $1$.

In order to rewrite the M-stationarity conditions of \eqref{eq:MPCC} as a Newton differentiable square system, 
we have to encode the complementarity-slackness-condition \eqref{eq:MPCC_stat_lambda} 
associated with the inequality constraints
as well as
the characteristic M-stationarity conditions \eqref{eq:MPCC_stat_W} and \eqref{eq:MPCC_stat_M}
with the aid of equations.
While the former can be achieved with the aid of so-called NCP-functions like the
minimum- or Fischer--Burmeister-function $\pi_\textup{min},\pi_\textup{FB}\colon\R^2\to\R$ given by
\begin{equation}\label{eq:NCP_functions}
	\pi_\textup{min}(a,b):=\min(a,b),\qquad
	\pi_\textup{FB}(a,b):=\sqrt{a^2+b^2}-a-b,
\end{equation}
see e.g.\ \cite{SunQi1999} for an introduction to and an overview of NCP-functions,
the latter is much more involved.
For some (yet unknown) function $\varphi\colon\R^4\to\R^2$, we consider the residual
$F\colon\R^n\times\R^r\times\R^s\times\R^t\times\R^t\to\R^n\times\R^r\times\R^s\times\R^{2t}$ given by
\begin{equation}\label{eq:min_residual}
	F(x,\lambda,\eta,\mu,\nu)
	:=
	\begin{bmatrix}
		\nabla_x L(x,\lambda,\eta,\mu,\nu)
		\\
		[\pi_\textup{min}(-g_i(x),\lambda_i)]_{I^r}
		\\
		h(x)
		\\
		[\varphi(G_i(x),H_i(x),\mu_i,\nu_i)]_{I^t}.
	\end{bmatrix}
\end{equation}
In order to obtain that $F(x,\lambda,\eta,\mu,\nu)=0$ holds if and only if 
$x$ is an M-stationary point of \eqref{eq:MPCC} such that $(\lambda,\eta,\mu,\nu)\in\Lambda_\textup{M}(x)$ is valid,
$\varphi$ has to satisfy
\[
	\varphi(a,b,\mu,\nu)=0
	\quad\Leftrightarrow\quad
	(a,b,\mu,\nu)\in M_\textup{M},
\]
where $M_\textup{M}$ has been defined in \eqref{eq:M_stat_set}.
Let us define $\psi_1,\psi_2,\psi_3,\varphi_1\colon\R^4\to\R$ by means of
\[
	\begin{aligned}
	\psi_1(a,b,\mu,\nu)&:=\max(-a,|b|,|\mu|),&
	\quad
	\psi_2(a,b,\mu,\nu)&:=\max(-b,|a|,|\nu|),&
	\\
	\psi_3(a,b,\mu,\nu)&:=\max(|a|,|b|,\mu,\nu),&
	\quad
	\varphi_1(a,b,\mu,\nu)&:=\min\limits_{j\in\{1,2,3\}}\psi_j(a,b,\mu,\nu).&
	\end{aligned}
\]
Noting that we have
\begin{align*}
	M_\textup{M}
	&=
	\set{(a,0,0,\nu)\,|\,a\geq 0}
	\cup
	\set{(0,b,\mu,0)\,|\,b\geq 0}
	\cup
	\set{(0,0,\mu,\nu)\,|\,\mu,\nu\leq 0},
\end{align*}
$M_\textup{M}$ crumbles into three convex branches which correspond to the
zero sublevel sets of the functions $\psi_1$, $\psi_2$, and $\psi_3$, respectively.
Hence, we have
\[
	\varphi_1(a,b,\mu,\nu)=0
	\quad\Leftrightarrow\quad
	(a,b,\mu,\nu)\in M_\textup{M},
\] 
see \cite[Lemma~3.1]{HarderMehlitzWachsmuth2021}. 
Unfortunately, the function $\varphi$ appearing in \eqref{eq:min_residual} needs to possess a second
component which, in order to avoid issues regarding Newton differentiability and
invertibility of the resulting Newton derivative, has to be chosen with care.

A suitable priority and chain rule can be used in order to guarantee that
the Newton derivative of the function $\varphi_1$ 
only possesses values in the set 
$\set{\pm\mathtt e_j\in\R^4 \given j\in\set{1,\ldots,4}}$,
see \cite[Example~2.8, Lemma~2.11]{HarderMehlitzWachsmuth2021}. Hence, we can define
$\varphi_2\colon\R^4\to\R$ by means of
\[
	\varphi_2(a,b,\mu,\nu)
	:=
	\begin{cases}
		\min(|b|,|\nu|)	&	D\varphi_1(a,b,\mu,\nu)=\pm\mathtt e_1,
		\\
		\min(|a|,|\mu|)	&	D\varphi_1(a,b,\mu,\nu)=\pm\mathtt e_2,
		\\
		|b|				&	D\varphi_1(a,b,\mu,\nu)=\pm\mathtt e_3,
		\\
		|a|				&	D\varphi_1(a,b,\mu,\nu)=\pm\mathtt e_4.
	\end{cases}
\]
The subsequently stated lemma summarizes some elementary properties of the function $\varphi$
which possesses the components $\varphi_1$ and $\varphi_2$ from above, 
see \cite[Lemma~3.3]{HarderMehlitzWachsmuth2021}.

\begin{lemma}\label{lem:properties_of_varphi}
			We have $\varphi(a,b,\mu,\nu)=0$ if and only if $(a,b,\mu,\nu)\in M_\textup{M}$.
			Furthermore, $\varphi$ is Newton differentiable of order $\infty$ on $M_\textup{M}$,
			and it is not continuous in any open neighborhood of $M_\textup{M}$.
\end{lemma}

Based on the chain rule, we obtain that the residual $F$ from \eqref{eq:min_residual}, 
with the function $\varphi$ characterized above, 
is Newton differentiable on its domain with an explicitly computable Newton derivative, 
and whenever the second-order derivatives of $f$, $g$, $h$, $G$, and $H$ are locally Lipschitzian,
then the order of Newton differentiability is $1$.
Hence, we can apply the abstract Newton scheme from \eqref{eq:Newton_iteration} to $F$.
One can check that the aforementioned local convergence guarantees apply to a given M-stationary
point $\bar x\in\R^n$ of \eqref{eq:MPCC} with associated multiplier 
$(\bar\lambda,\bar\eta,\bar\mu,\bar\nu)\in\Lambda_\textup{M}(\bar x)$ whenever
MPCC-LICQ and MPCC-SSOC (with respect to $(\bar\lambda,\bar\eta,\bar\mu,\bar\nu)$) are valid at $\bar x$, 
see \cite[Theorems~4.1, 4.2]{HarderMehlitzWachsmuth2021}.
Let us recall that due to \cref{lem:properties_of_varphi}, the residual $F$ is discontinuous, 
so most of the standard nonsmooth Newton-type methods do not apply to the situation at hand.

\subsubsection{The algorithm}\label{sec:SSN_alg}

Without a suitable globalization strategy, a pure Newton-type scheme is often not competitive.
Typically, one uses gradient steps with respect to the squared Euclidean norm of the residual
in combination with a line search procedure whenever the Newton step is not well defined or
does not yield a sufficient decrease in the norm of the residual. In our setting, we want to
rely on a similar strategy. However, for the residual $F$ from \eqref{eq:min_residual}, the
associated function $\tfrac12\norm{F(\cdot)}^2$ is not continuously differentiable.
In order to circumvent this issue, we define yet another residual of the M-stationarity system
whose squared norm is smooth. Therefore, we heavily rely on the Fischer--Burmeister function
defined in \eqref{eq:NCP_functions} whose square is known to be continuously differentiable,
see e.g.\ \cite[Proposition~3.4]{FacchineiSoares1997}.

Let us define $\theta_1,\theta_2,\theta_3,\theta_4\colon\R^4\to\R$ 
by means of
\[
	\begin{aligned}
		\theta_1(a,b,\mu,\nu)&:=|\pi_\textup{FB}(a,b)|,&
		\quad
		\theta_2(a,b,\mu,\nu)&:=\pi_\textup{FB}(|a|,|\mu|),&
		\\
		\theta_3(a,b,\mu,\nu)&:=\pi_\textup{FB}(|b|,|\nu|),&
		\quad
		\theta_4(a,b,\mu,\nu)&:=
			\begin{cases}
				0								&	\mu,\nu\leq 0 \\ 
				\pi_\textup{FB}(|\mu|,|\nu|),	&	\text{else,} 
			\end{cases}
			&
	\end{aligned}
\]
and let $\theta\colon\R^4\to\R^4$ be the function which possesses the component functions $\theta_1,\theta_2,\theta_3,\theta_4$.
Using the aforementioned properties of the Fischer--Burmeister function, it is not hard to see
that $\tfrac12\norm{\theta(\cdot)}^2$ is continuously differentiable.
Furthermore, one can check that 
\[
	\theta(a,b,\mu,\nu)=0\quad\Leftrightarrow\quad(a,b,\mu,\nu)\in M_\textup{M}.
\] 
Next, define $F_\textup{FB}\colon\R^n\times\R^r\times\R^s\times\R^t\times\R^t\to\R^n\times\R^r\times\R^s\times\R^{4t}$ by means of
\[
	F_\textup{FB}(x,\lambda,\eta,\mu,\nu)
	:=
	\begin{bmatrix}
		\nabla_x L(x,\lambda,\eta,\mu,\nu)
		\\
		[\pi_\textup{FB}(-g_i(x),\lambda_i)]_{I^r}
		\\
		h(x)
		\\
		[\theta(G_i(x),H_i(x),\mu_i,\nu_i)]_{I^t}
	\end{bmatrix}
	,
\]
and observe that $F_\textup{FB}(x,\lambda,\eta,\mu,\nu)=0$ holds if and only if 
$x\in\R^n$ is an M-stationary point of \eqref{eq:MPCC} with $(\lambda,\eta,\mu,\nu)\in\Lambda_\textup{M}(x)$.
Additionally, with the aid of \cite[Lemma~3.1]{Tseng1996}, it is possible to show
the existence of constants $c,C>0$ such that
\begin{equation}\label{eq:equivalence_of_residuals}
	c\norm{F_\textup{FB}(x,\lambda,\eta,\mu,\nu)}^2
	\leq
	\norm{F(x,\lambda,\eta,\mu,\nu)}^2
	\leq
	C\norm{F_\textup{FB}(x,\lambda,\eta,\mu,\nu)}^2,
\end{equation}
i.e., the residuals $F$ and $F_\textup{FB}$ are, to some extent, equivalent, see \cite[Lemma~5.1]{HarderMehlitzWachsmuth2021}.

By construction, the map $\Phi_\textup{FB}\colon\R^n\times\R^r\times\R^s\times\R^t\times\R^t\to\R$ given by
\[
	\Phi_\textup{FB}(x,\lambda,\eta,\mu,\nu)
	:=
	\tfrac12\norm{F_\textup{FB}(x,\lambda,\eta,\mu,\nu)}^2
\]
is now continuously differentiable, and we will use it for the globalization of our Newton method.
To keep the notation short, we introduce a surrogate variable
\[
	\bz:=(x,\lambda,\eta,\mu,\nu)\in\R^n\times\R^r\times\R^s\times\R^t\times\R^t,
\]
and similarly, $\bz^k$ for some iteration index $k\in\N$ has to be be understood.
Based on the function $\Phi_\textup{FB}$ from above and a given iterate $\bz^k$, 
our globalization strategy works as follows, see \cite[Section~3]{DeLucaFacchineiKanzow2000}. 
We first check whether the Newton direction $\bd^k$ can be computed as a solution of the Newton system
\begin{equation}\label{eq:Newton_system}
	DF(\bz^k)\bd^k=-F(\bz^k)
\end{equation}
and satisfies the ratio test
\begin{equation}\label{eq:ratio_test}
	\Phi_\textup{FB}(\bz^k+\bd^k)\leq q_\textup{nsn}\,\Phi_\textup{FB}(\bz^k)
\end{equation}
for some given $q_\textup{nsn}\in(0,1)$. In this case, we do a full Newton step $\bz^{k+1}:=\bz^k+\bd^k$.
In case where the Newton system \eqref{eq:Newton_system} cannot be solved or if its solution 
does not pass a standard angle test, we simply set $\bd^k:=-\nabla\Phi_\textup{FB}(\bz^k)$.
Afterwards, we use an Armijo line search to obtain a step size $\alpha_k>0$ and set $\bz^{k+1}:=\bz^k+\alpha_k\bd^k$.

In \cref{alg:SSN}, the pseudocode of our nonsmooth Newton method is stated.
We would like to mention that \cref{alg:SSN} is a descent method with respect to $\Phi_\textup{FB}$,
and that we would like to compute the global minimizers of this function as these points coincide 
with the solution set of the M-stationarity conditions of \eqref{eq:MPCC}.
Due to \eqref{eq:equivalence_of_residuals}, it is also clear that $\Phi_\textup{FB}(\bz^k)\to 0$
directly yields $\norm{F(\bz^k)}\to 0$, i.e., the termination criterion in \cref{alg:SSN} is reasonable.
In a practical implementation of \cref{alg:SSN}, one should also check if $\norm{\nabla\Phi_\textup{FB}(\bz^k)}$ 
becomes small since accumulation points of sequences generated by descent methods often turn out to be merely
stationary points of the underlying function, and a similar observation can be made for \cref{alg:SSN}.

\begin{algorithm2e}[htb]
	\SetAlgoLined
	\KwData{parameters $q_\textup{nsn},\tau_{\textup{nsn}},\rho, \sigma, \beta \in(0,1)$,
	starting point $\bz^0\in\R^n\times\R^r\times\R^s\times\R^t\times\R^t$}
	Set $k:=0$\;
	\While{$\norm{F(\bz^k)}>\tau_{\textup{nsn}}$}{
		Solve \eqref{eq:Newton_system}\;
		\eIf{$\bd^k$ is well defined \KwAnd ratio test \eqref{eq:ratio_test} is satisfied}{
			Set $\bz^{k+1}:=\bz^k+\bd^k$\;
		}{
			\If{$\bd^k$ is not well defined \KwOr $\nabla\Phi_\textup{FB}(\bz^k)^\top \bd^k > -\rho\norm{\bd^k}\norm{\nabla\Phi_\textup{FB}(\bz^k)}$}{
				Set $\bd^k := -\nabla\Phi_\textup{FB}(\bz^k)$\;
			}
			Determine $\bz^{k+1}:=\bz^k+\alpha_k\bd^k$ using an Armijo line search for $\Phi_\textup{FB}$, 
			i.e., $\alpha_k := \beta^{i_k}$, where $i_k \in \N_0$ is the smallest non-negative integer with
			$\Phi_\textup{FB}(\bz^k + \beta^{i_k} \, \bd^k) \le \Phi_\textup{FB}(\bz^k) + \sigma \, \beta^{i_k} \, \nabla\Phi_\textup{FB}(\bz^k)^\top \bd^k$\;
		}
		Set $k := k+1$\;
	}
	\caption{Globalized nonsmooth Newton method for \eqref{eq:MPCC}.}
	\label{alg:SSN}
\end{algorithm2e}

\subsubsection{Convergence guarantees}\label{sec:SSN_convergence}

Let us present some facts about the convergence behavior of \cref{alg:SSN}, 
taken from \cite[Theorem~5.2]{HarderMehlitzWachsmuth2021}.
Therefore, we assume that \cref{alg:SSN} produces an infinite sequence $\{\bz^k\}_{k\in\N}$.
It is clear that, whenever \eqref{eq:ratio_test} is satisfied infinitely often, then
$\Phi_\textup{FB}(\bz^k)\to 0$ must be valid as \cref{alg:SSN} is a descent method.
Consequently, in this situation, each accumulation point of $\{\bz^k\}_{k\in\N}$ is an
M-stationary point of \eqref{eq:MPCC} together with an associated multiplier in this situation.
In more general situations, one can merely verify that each accumulation point of $\{\bz^k\}_{k\in\N}$
is a stationary point of $\Phi_\textup{FB}$. 
However, if the primal component of some accumulation point of the sequence satisfies MPCC-LICQ and MPCC-SSOC 
with respect to the associated dual components, then the overall sequence converges superlinearly
to this point since, for all sufficiently large iterations, the full Newton step is accepted,
and \cref{alg:SSN} behaves like a local nonsmooth Newton method. 
Whenever the second-order derivatives of all data functions in \eqref{eq:MPCC}
are locally Lipschitzian, then the convergence of the whole sequence is already quadratic.

Finally, let us mention that in situations where $f$ is quadratic while $g$, $h$, $G$, and $H$ are affine,
it is possible to show local fast convergence of \cref{alg:SSN} under MPCC-SSOC and a slightly weaker
constraint qualification than MPCC-LICQ. 
More precisely, one only needs the linear independence of the gradients
\[
	\nabla g_i(\bar x)\,(i\in I_g^+),
	\;
	\nabla h_i(\bar x)\,(i\in I^s),
	\;
	\nabla G_i(\bar x)\,(i\in I_{0+}\cup I_{00}^{\pm\R}),
	\;
	\nabla H_i(\bar x)\,(i\in I_{+0}\cup I_{00}^{\R\pm})
\]
in this situation. However, as the appearing index sets depend on the associated multipliers 
which solve the system of M-stationarity, this condition is more difficult to check.
Detailed comments can be found in \cite[Section~6]{HarderMehlitzWachsmuth2021}.

\section{Numerical results for the inverse optimal control problem}\label{sec:numerics_IOC}

For our numerical examination of the methods presented in \cref{sec:computing_M_stationary_points},
we aim to challenge both algorithms by an instance of the inverse optimal control
problem \eqref{eq:IOC} discussed in \cite[Section~5.4]{HarderWachsmuth2022}.
We choose $\Omega:=(0,1)^2$, $Y:=\R$, $y_\textup{d}:=0$, $\alpha:=10^{-3}$, $u_\textup{a}\equiv 0$, 
$\zeta\equiv 1$, and $c,S\colon L^2(\Omega)\to\R$ are given by
\[
	\forall u\in L^2(\Omega)\colon\quad
	c(u) := \frac12\norm{u-u_\textup{o}}^2_{L^2(\Omega)},
	\qquad
	Su := \innerprod{1}{u}_{L^2(\Omega)},
\] 
where $u_\textup{o}\equiv 1$ is used for the observed control.
The upper-level lower bound $w_\textup{a}$ will vary throughout the experiments.
We note that, due to \cref{thm:mstat}, all local minimizers of this problem
satisfy the pointwise M-stationarity conditions of the associated MPCC
\eqref{eq:KKTR}.

We aim to solve \eqref{eq:IOC} based on its KKT reformulation \eqref{eq:KKTR}.
The latter problem is discretized by a standard finite element approach.
More precisely, we choose a suitable triangulation of $\Omega$ into $128$ triangles, 
and the variables $u$ and $\xi$ are discretized by piecewise constant functions while
the function $w$ is discretized with the aid of piecewise linear and continuous functions.
The discretized complementarity-constrained optimization problem
associated with \eqref{eq:KKTR} is a linear-quadratic MPCC and, thus, its local
minimizers are M-stationary as MPCC-GCQ is valid at each feasible point.

\cref{alg:ALM,alg:SSN} have been implemented in MATLAB2022b.
In \cref{alg:ALM}, the appearing parameters are set to
$\gamma:=10$,
$q_\textup{alm}:=0.8$,
$C:=10^{20}$,
$\varepsilon_k:=10^{-4}\,(k+1)^{-1/2}$,
and
$\tau_\textup{alm}:=10^{-5}$.
The initial penalty parameter $\rho_0$ is computed as described in \cite[Section~6]{JiaKanzowMehlitzWachsmuth2023}.
The computational solution of the appearing subproblems is realized 
via the (nonmonotone) projected gradient method from \cite[Section~3]{JiaKanzowMehlitzWachsmuth2023}
with parameters chosen according to \cite[Section~6]{JiaKanzowMehlitzWachsmuth2023}.
As the discretized MPCC under consideration comes along with simple enough structure,
we implemented \cref{alg:ALM} in the flavor where slack variables are avoided and projections
onto the associated set $\mathcal D$ from \eqref{eq:definition_of_D} are utilized.
In our experiments, we made use of
$q_\textup{nsn}:=0.999$,
$\tau_\textup{nsn}:=10^{-11}$,
$\rho:=10^{-3}$,
and
$\sigma:=\beta:=0.5$ in \cref{alg:SSN}.
The maximum number of iterations for both algorithms is set to $1000$.
Starting points for both algorithms are constructed as follows:
the actual variables are initialized as random vectors with entries chosen via a standard normal distribution,
the initial guess for the multipliers is the all-zero vector.

For our first experiment, we choose $w_\textup{a}\equiv 0$.
It has been shown in \cite[Section~5.4]{HarderWachsmuth2022} that the
associated uniquely determined global minimizer $(\bar u,\bar\xi,\bar w):=(0,0,0)$
is pointwise M- but not S-stationary.
The associated (upper-level) objective function value of the original and discretized
problem is $2$.
We ran both algorithms based on (the same) $10$ random starting points.
The results are presented in \cref{tab:Ex1} where we list the number of iterations,
the (upper-level) function value of the final iterate, the computation time in seconds,
the final value of the penalty parameter for \cref{alg:ALM},
and the number of full Newton steps as well as gradient steps done by \cref{alg:SSN}.

\begin{table}[ht]
\centering
\begin{tabular}{cccccccccc}
\toprule
&\multicolumn{4}{c}{\cref{alg:ALM}}&\multicolumn{5}{c}{\cref{alg:SSN}}\\
\cmidrule(r){2-5} \cmidrule{6-10}
& \# iterations & value & time & $\rho$ & \# iterations & value & time & \# full steps & \# gradient steps\\
\midrule
1 & 15 & 2.00 & 9.78 & 1.59$\cdot 10^8$ & 11 & 2.00 & 0.55 & 6 & 5\\
2 & 13 & 2.00 & 3.76 & 9.19$\cdot 10^6$ & 7 & 2.00 & 0.29 & 4 & 3\\
3 & 18 & 2.00 & 7.77 & 6.99$\cdot 10^6$ & 5 & 2.00 & 0.17 & 3 & 2\\
4 & 16 & 2.00 & 5.22 & 7.11$\cdot 10^6$ & 5 & 2.00 & 0.16 & 3 & 2\\
5 & 13 & 2.00 & 5.12 & 4.44$\cdot 10^6$ & 5 & 2.00 & 0.17 & 3 & 2\\
6 & 16 & 2.00 & 6.87 & 7.23$\cdot 10^6$ & 8 & 2.00 & 0.29 & 4 & 4\\
7 & 14 & 2.00 & 3.93 & 7.85$\cdot 10^6$ & 6 & 2.00 & 0.20 & 4 & 2\\
8 & 15 & 2.00 & 10.00 & 7.52$\cdot 10^7$ & 8 & 2.00 & 0.28 & 5 & 3\\
9 & 16 & 2.00 & 6.12 & 7.29$\cdot 10^6$ & 7 & 2.00 & 0.26 & 4 & 3\\
10 & 16 & 2.00 & 7.08 & 7.94$\cdot 10^7$ & 5 & 2.00 & 0.18 & 3 & 2\\
\bottomrule
\end{tabular}
\caption{Results of experiments for $w_\textup{a}\equiv 0$.}
\label{tab:Ex1}
\end{table}

We immediately see that both algorithms find the global minimizer in all runs.
On the one hand, \cref{alg:ALM} needs between $13$ and $16$ iterations until the termination criterion is hit,
and this takes between $3$ and $10$ seconds of time.
During each run, the penalty parameter is enlarged $6$ to $8$ times.
On the other hand, \cref{alg:SSN} terminates after at most $11$ iterations out of which the most ones
are full Newton steps. Less than $0.6$ seconds of time are needed in each run.
Hence, we guess that the first steps of \cref{alg:SSN} drive the iterates into the radius
of attraction associated with the global minimizer (together with suitable multipliers)
of Newton's method,
so that local fast convergence can be observed during the last steps.
In this regard, \cref{alg:SSN} outruns \cref{alg:ALM} for this experiment.
Let us note that the running time of \cref{alg:ALM} is dominated by the running time of the
subproblem solver. As the latter one is a projected gradient method, 
it naturally runs a lot of (very cheap) iterations.

For our second experiment, we choose $w_\textup{a}\equiv-0.05$.
As this enlarges the feasible set, the (upper-level) objective function value 
of the associated global minimizer is not larger than $2$.
Again, we ran both algorithms for the same $10$ random starting points 
used before and obtained the results stated in \cref{tab:Ex2}.

\begin{table}[ht]
\centering
\begin{tabular}{cccccccccc}
\toprule
&\multicolumn{4}{c}{\cref{alg:ALM}}&\multicolumn{5}{c}{\cref{alg:SSN}}\\
\cmidrule(r){2-5} \cmidrule{6-10}
& \# iterations & value & time & $\rho$ & \# iterations & value & time & \# full steps & \# gradient steps\\
\midrule
1 & 14 & 1.88 & 595.08 & 1.74$\cdot 10^8$ & 1000 & 1.87 & 11.85 & 1 & 941\\
2 & 13 & 1.88 & 102.41 & 1.00$\cdot 10^8$ & 1000 & 1.87 & 11.12 & 1 & 963\\
3 & 14 & 1.88 & 117.34 & 7.48$\cdot 10^7$ & 1000 & 1.87 & 11.42 & 1 & 932\\
4 & 14 & 1.88 & 93.26 & 7.68$\cdot 10^7$ & 1000 & 1.87 & 12.42 & 1 & 817\\
5 & 14 & 1.88 & 82.99 & 4.78$\cdot 10^7$ & 1000 & 1.87 & 10.87 & 1 & 974\\
6 & 14 & 1.88 & 72.39 & 7.86$\cdot 10^7$ & 1000 & 1.87 & 13.31 & 1 & 970\\
7 & 16 & 1.88 & 358.47 & 8.38$\cdot 10^7$ & 1000 & 1.87 & 11.96 & 1 & 820\\
8 & 13 & 1.88 & 66.18 & 8.02$\cdot 10^7$ & 1000 & 1.87 & 10.73 & 1 & 871\\
9 & 14 & 1.88 & 66.99 & 7.78$\cdot 10^7$ & 1000 & 1.86 & 10.49 & 1 & 932\\
10 & 14 & 1.88 & 115.14 & 8.66$\cdot 10^7$ & 1000 & 1.87 & 10.72 & 1 & 940\\
\bottomrule
\end{tabular}
\caption{Results of experiments for $w_\textup{a}\equiv -0.05$.}
\label{tab:Ex2}
\end{table}

While the behavior of \cref{alg:ALM} does not change significantly
with respect to iteration numbers and the evolution of the penalty parameter,
computation time drastically increases.
\Cref{alg:ALM} now needs between $1$ and $10$ minutes to terminate
which, as already mentioned, is due to large iteration numbers of the subproblem solver.
\Cref{alg:SSN} terminates in each of the runs since the maximum number of iterations is reached.
Just one full Newton step is carried out in each run, 
and most of the remainder steps are gradient steps. 
This slows down the convergence so that the underlying termination
criterion is not hit within the maximum number of iterations.
However, based on the computed function value,
we guess that the final iterate produced by both methods
in each of the runs is close to the global minimizer 
of the (perturbed) inverse optimal control problem.

In our third experiment, we aim to combine \cref{alg:ALM,alg:SSN}
in order to obtain a method which benefits, on the one hand, from 
the robust behavior of \cref{alg:ALM} and, on the other hand,
from the potential local fast convergence of \cref{alg:SSN}.
Therefore, for $w_\textup{a}\equiv -0.05$ again,
we first run \cref{alg:ALM} on the set of $10$ random starting
points already used before, but with $\tau_\textup{alm} := 10^{-5}$, i.e.,
with a less stringent termination criterion.
Second, the final iterates are then used as starting points for
\cref{alg:SSN}. \cref{tab:Ex2_warmstart} clearly documents that
this warm starting is promising as it reduces the overall
number of iterations and the overall computation time since
\cref{alg:SSN} terminates after just one full Newton step in each
of the runs. 
	In the light of \cite[Theorem~2.9]{HarderMehlitzWachsmuth2021},
	the latter behavior is not surprising as we aim to solve
	a linear-quadratic complementarity-constrained problem
	while \cref{alg:ALM} seemingly drives the sequence of
	iterates in the radius of local fast convergence of \cref{alg:SSN}.

\begin{table}[ht]
\centering
\begin{tabular}{cccccccccc}
\toprule
&\multicolumn{4}{c}{\cref{alg:ALM}}&\multicolumn{5}{c}{\cref{alg:SSN}}\\
\cmidrule(r){2-5} \cmidrule{6-10}
 & \# iterations & value & time & $\rho$ & \# iterations & value & time & \# full steps & \# gradient steps\\
\midrule
1 & 10 & 1.88 & 120.37 & 1.74$\cdot 10^6$ & 1 & 1.88 & 0.02 & 1 & 0\\
2 & 11 & 1.88 & 75.55 & 1.00$\cdot 10^7$ & 1 & 1.88 & 0.01 & 1 & 0\\
3 & 12 & 1.88 & 76.75 & 7.48$\cdot 10^6$ & 1 & 1.88 & 0.01 & 1 & 0\\
4 & 12 & 1.88 & 65.02 & 7.68$\cdot 10^6$ & 1 & 1.88 & 0.01 & 1 & 0\\
5 & 10 & 1.88 & 40.64 & 4.78$\cdot 10^5$ & 1 & 1.88 & 0.01 & 1 & 0\\
6 & 12 & 1.88 & 42.74 & 7.86$\cdot 10^6$ & 1 & 1.88 & 0.01 & 1 & 0\\
7 & 12 & 1.88 & 101.89 & 8.38$\cdot 10^5$ & 1 & 1.88 & 0.01 & 1 & 0\\
8 & 11 & 1.88 & 37.17 & 8.02$\cdot 10^6$ & 1 & 1.88 & 0.01 & 1 & 0\\
9 & 12 & 1.88 & 41.36 & 7.78$\cdot 10^6$ & 1 & 1.88 & 0.01 & 1 & 0\\
10 & 12 & 1.88 & 68.78 & 8.66$\cdot 10^6$ & 1 & 1.88 & 0.01 & 1 & 0\\
\bottomrule
\end{tabular}
\caption{Results of experiments for $w_\textup{a}\equiv -0.05$ via warm starting.}
\label{tab:Ex2_warmstart}
\end{table}

	Our final experiment aims to visualize the numerical behavior of the
	warm started solution method on a finer mesh.
	Therefore, we reconsider the problem of interest with $w_\textup{a}\equiv-0.05$
	on a refined grid of 512 triangles.
	Again, we run the algorithm on a (due to the refinement of the mesh)
	different set of 10 random starting points, 
	and make use of $\tau_\textup{alm}:=10^{-6}$.
	Additionally, we immediately abort \cref{alg:ALM} if a maximum number
	of $10^7$ cumulated function evaluations of the augmented Lagrangian function is
	exceeded but still hand the last iterate over to \cref{alg:SSN}.
	This happened for the starting points indexed by 2, 4, 5, 6, and 8.
	The results are documented in \cref{tab:Ex2_warmstart_refined}.	
	On the one hand, due to the larger size of the problem, we observe higher 
	computational times for \cref{alg:ALM,alg:SSN}. 
	On the other hand, the total number iterations for \cref{alg:ALM}
	and the final value of the penalty parameter increase just slightly,
	and \cref{alg:SSN} still shows one-step convergence. 
	
	\begin{table}[ht]
	\centering
	\begin{tabular}{cccccccccc}
	\toprule
	&\multicolumn{4}{c}{\cref{alg:ALM}}&\multicolumn{5}{c}{\cref{alg:SSN}}\\
	\cmidrule(r){2-5} \cmidrule{6-10}
 	& \# iterations & value & time & $\rho$ & \# iterations & value & time & \# full steps & \# gradient steps\\
	\midrule
	1 & 16 & 1.88 & 446.57 & 1.04$\cdot 10^8$ & 1 & 1.88 & 0.15 & 1 & 0\\
	2 & 16 & 1.88 & 2488.56 & 9.00$\cdot 10^7$ & 1 & 1.88 & 0.13 & 1 & 0\\
	3 & 19 & 1.88 & 694.89 & 8.03$\cdot 10^7$ & 1 & 1.88 & 0.13 & 1 & 0\\
	4 & 18 & 1.88 & 2924.33 & 7.51$\cdot 10^7$ & 1 & 1.88 & 0.12 & 1 & 0\\
	5 & 17 & 1.88 & 3052.80 & 6.90$\cdot 10^7$ & 1 & 1.88 & 0.12 & 1 & 0\\
	6 & 16 & 1.88 & 2707.00 & 8.72$\cdot 10^7$ & 1 & 1.88 & 0.12 & 1 & 0\\
	7 & 16 & 1.88 & 909.23 & 8.78$\cdot 10^7$ & 1 & 1.88 & 0.12 & 1 & 0\\
	8 & 19 & 1.88 & 3565.73 & 8.50$\cdot 10^7$ & 1 & 1.88 & 0.12 & 1 & 0\\
	9 & 18 & 1.88 & 1080.83 & 6.99$\cdot 10^8$ & 1 & 1.88 & 0.13 & 1 & 0\\
	10 & 15 & 1.88 & 984.20 & 7.29$\cdot 10^7$ & 1 & 1.88 & 0.12 & 1 & 0\\
	\bottomrule
	\end{tabular}
	\caption{Results of experiments for $w_\textup{a}\equiv -0.05$ via warm starting on a refined mesh.}
	\label{tab:Ex2_warmstart_refined}
	\end{table}

\section{Conclusions}\label{sec:conclusions}

We demonstrated by means of a class of inverse optimal control problems 
that pointwise M-stationarity may yield a reasonable necessary optimality condition
for MPCCs in Lebesgue spaces
in situations where pointwise S-stationarity fails.
For this purpose, we relied on geometric arguments inspired by some recent progress
in the field of finite-dimensional complementarity-constrained optimization.
It is an interesting question of future research whether the presented
techniques can be refined in order to obtain pointwise M-stationarity
for a broader class of MPCCs in Lebesgue spaces, and if it is possible to extend
the approach to MPCCs in more difficult function spaces like Sobolev spaces.

Furthermore, we reviewed two recently developed algorithms for the computation
of M-stationary points of finite-dimensional MPCCs - an augmented Lagrangian and a Newton-type method.
We challenged both algorithms by (discretized) instances 
of the aforementioned inverse optimal control problem. 
While the augmented Lagrangian method turned out to be very robust against perturbations of the starting
point, we were in position to observe local fast convergence of the Newton-type method in some situations.
It has been demonstrated that using the augmented Lagrangian method for warm starting of
the Newton method is beneficial with respect to iteration numbers and computation time.
It remains to be seen whether the ideas used for the construction of these algorithms
can be extended to the infinite-dimensional setting.

%\bibliographystyle{spmpsci}
%\bibliography{references}

\end{document}